\documentclass[twocolumn]{autart}    

\usepackage{graphicx}          

\usepackage{epsfig} 
\usepackage{amsfonts}
\usepackage{amsmath,amssymb}
\usepackage{comment}
\usepackage{color}
\usepackage{cases}
\usepackage[subrefformat=parens]{subcaption}

\usepackage{standalone}
\usepackage{algpseudocode}
\usepackage{hyperref}
\usepackage{algorithm}
\usepackage{array}
\usepackage{textcomp}
\usepackage{stfloats}
\usepackage{url}
\usepackage{verbatim}
\usepackage{graphicx}

\usepackage{cite}

\usepackage{mathtools}				 
\usepackage{enumitem}

\usepackage{tikz}

\usepackage{bm}
\usepackage{here}
\usepackage{mathrsfs}

\usepackage{stmaryrd}

\makeatletter
\newcommand{\specialcell}[1]{\ifmeasuring@#1\else\omit$\displaystyle#1$\ignorespaces\fi}
\makeatother

\makeatletter
\newcommand{\pushright}[1]{\ifmeasuring@#1\else\omit\hfill$\displaystyle#1$\fi\ignorespaces}
\newcommand{\pushleft}[1]{\ifmeasuring@#1\else\omit$\displaystyle#1$\hfill\fi\ignorespaces}
\makeatother

\newtheorem{problem}{Problem}
\newtheorem{proposition}{Proposition}

\newtheorem{lemma}{Lemma}
\newtheorem{theorem}{Theorem}
\newtheorem{corollary}{Corollary}
\newtheorem{rmk}{Remark}

\newtheorem{assumption}{Assumption}

\newcommand{\calE}{{\mathcal E}}    
    
\newcommand{\calG}{{\mathcal G}}    
\newcommand{\calH}{{\mathcal H}}

\newcommand{\calM}{{\mathcal M}}

\newcommand{\calQ}{{\mathcal Q}}    
\newcommand{\calR}{{\mathcal R}}    
    
\newcommand{\calT}{{\mathcal T}}

\newcommand{\calX}{{\mathcal X}}

\newcommand{\scrG}{{\mathscr G}}

\newcommand{\bbR}{{\mathbb R}}

\newcommand{\bbU}{{\mathbb U}}

\newcommand{\bbX}{{\mathbb X}}    
    
\newcommand{\bbZ}{{\mathbb Z}}

\newcommand{\argmin}{\mathop{\rm arg~min}\limits}

\newcommand{\bfa}{{\mathbf a}}
\newcommand{\bfb}{{\mathbf b}}

\newcommand{\rmd}{{\rm d}}

\newcommand{\rmh}{{\rm h}}

\newcommand{\rmH}{{\rm H}}

\newcommand{\sfd}{{\sf d}}
\newcommand{\sfe}{{\sf e}}

\newcommand{\bff}{{\bf f}}
\newcommand{\bfg}{{\bf g}}
\newcommand{\bfT}{{\bf T}}
\newcommand{\bfone}{{\bf 1}}

\newcommand{\bmA}{{\boldsymbol A}}
\newcommand{\bmB}{{\boldsymbol B}}
\newcommand{\bmC}{{\boldsymbol C}}
\newcommand{\bmG}{{\boldsymbol G}}
\newcommand{\bmK}{{\boldsymbol K}}
\newcommand{\bmP}{{\boldsymbol P}}

\newcommand{\bmc}{{\boldsymbol c}}
\newcommand{\bmell}{{\boldsymbol \ell}}
\newcommand{\bmu}{{\boldsymbol u}}
\newcommand{\bmx}{{\boldsymbol x}}

\newcommand{\Ee}{{\rm e}}   

\newcommand{\bbra}[1]{\ensuremath{[\![#1]\!]} }  

\newcommand{\dhil}{{d_\calH}}

\newcommand{\diagbox}{\boxbslash}

\newcommand{\one}{\scalebox{0.6}{\bf 1}_{\scalebox{0.4}{$N$}}}

\newcommand{\mpc}{{\rm MPC}}
\newcommand{\tmp}{{\rm tmp}}

\newcommand{\magenta}[1]{{\color{black}{#1}}} 
\newcommand{\blue}[1]{{\color{black}{#1}}} 
\newcommand{\rev}[1]{{\color{black}{#1}}} 
\newcommand{\fin}[1]{{\color{black}{#1}}} 
\newcommand{\red}[1]{{\color{black}{#1}}} 
\newcommand{\modi}[1]{{\color{black}{#1}}} 
\newcommand{\comm}[1]{{\color{black}{#1}}} 
\newcommand{\subm}[1]{{\color{black}{#1}}} 
\newcommand{\res}[1]{{\color{black}{#1}}} 
\newcommand{\blu}[1]{{\color{black}{#1}}} 
\newcommand{\las}[1]{{\color{black}{#1}}} 

\begin{document}

\begin{frontmatter}

\title{Entropic Model Predictive Optimal Transport\\over Dynamical Systems\thanksref{footnoteinfo}} 

\thanks[footnoteinfo]{This paper was not presented at any IFAC
	meeting. Corresponding author K.~Kashima. Tel. +81-75-753-5512.}

\author[Tokyo]{Kaito Ito}\ead{ka.ito@c.titech.ac.jp},    
\author[Kyoto]{Kenji Kashima}\ead{kk@i.kyoto-u.ac.jp}  

\address[Tokyo]{School of Computing, Tokyo Institute of Technology, Yokohama 226-8502, Japan}             
\address[Kyoto]{Graduate School of Informatics, Kyoto University, Kyoto 606-8501, Japan}  


\begin{abstract}                   
	We consider the optimal control problem of steering an agent population to a desired distribution over an infinite horizon.
  This is an optimal transport problem over dynamical \blu{systems}, which is challenging due to its high computational cost.
  In this paper, by using entropy regularization, we propose {\em Sinkhorn MPC}, which is a dynamical transport algorithm integrating model predictive control \subm{(MPC)} and the so-called Sinkhorn algorithm.
  The notable feature of the proposed method is that it achieves cost-effective transport in real time by performing control and transport planning simultaneously, which is illustrated in numerical examples.
	Moreover, under some assumption on iterations of the Sinkhorn algorithm integrated in MPC, we reveal the global convergence property for Sinkhorn MPC thanks to the entropy regularization.
  Furthermore, focusing on \fin{a quadratic control cost}, without the aforementioned assumption we show the ultimate boundedness and the local asymptotic stability for Sinkhorn MPC.

\end{abstract}

\end{frontmatter}

\section{Introduction}\label{sec:intro}
The problem of controlling a large number of agents has become a more and more important area in control theory with a view to applications in sensor networks, smart grids, intelligent transportation systems, and systems biology, to name a few~\cite{Chung2018,Chowdhury2015,Alasseur2020}.
One of the most fundamental tasks in this problem is to stabilize a collection of agents to a desired distribution shape with minimum cost.
This can be formulated as an optimal transport (OT) problem~\cite{Villani2003} between the empirical distribution based on the state of the agents and the target distribution over dynamical \blu{systems}.
The OT problem over dynamical systems consists of finding an assignment of agents to targets and control inputs that drive the agents to the assigned targets in order to minimize the total cost of interest.
The difficulty of this problem lies in the large scale of the collective dynamics.

{\it Literature review:}
The assignment problem has been extensively studied in the context of combinatorial optimization, and many methods to find the optimal assignment have been proposed such as the well-known Hungarian algorithm~\cite{Kuhn1955} and auction algorithm~\cite{Bertsekas1992}.
These algorithms have been applied to multi-agent assignment problems; see e.g., \cite{Yu2014,Mosteo2017} and references therein.
In the literature, the dynamics of agents are simplified as the single integrator dynamics, and easily computable assignment costs, e.g., distance-based cost, are considered in general.
On the other hand, when considering more general dynamics and cost functions for the stabilization to targets, it is difficult to obtain the associated assignment costs and optimal controls. This is because, in most cases, infinite horizon optimal control (OC) problems stabilizing agents to desired targets are computationally intractable.
A promising approach to overcome this problem is model predictive control~(MPC)~\cite{Mayne2014}, in which the current control input is determined by solving, at each sampling instant, a finite horizon OC problem using the current state as the initial state. For example \subm{in \cite{Morgan2016},} MPC is used to solve a finite horizon assignment problem over \subm{dynamical systems} in real time.
Now it is important to emphasize that when performing MPC for \subm{a dynamic} OT problem, \subm{it is desirable to update the target assignment for agents at each time as well as control inputs.}
However, when the number of the agents is large, solving the assignment problem at each sampling instant is computationally very expensive even with \subm{the Hungarian algorithm}.
\subm{Even worse, the changes of the assignment along the controlled state trajectories are not continuous, and this makes it difficult to ensure the stability of the dynamics under MPC.}

On the other hand, recently, a different approach to solve a dynamical assignment problem using OT theory has attracted much attention~\cite{Chen2018,Bakshi2020,Krishnan2018}.
In this approach, a large population limit is considered, and infinitely many agents are represented as a probability density of the state of a single system. 
Then, the dynamical assignment problem boils down to a density control problem~\cite{Chen2017,de2021,Ito2022bridge} finding a feedback control law that steers an initial \subm{state} density to a target density with minimum cost.
Consequently, this approach can avoid the difficulty due to the large scale of the collective dynamics.
\blu{Nevertheless}, it has the drawback that even for linear systems, the density control requires to solve a nonlinear partial differential equation such as the Monge-Amp\`{e}re equation or the Hamilton-Jacobi-Bellman equation, which are generally difficult to solve.

{\it Contributions:}
With this in mind, we deal with the collective dynamics directly without taking the number of agents to infinity, but utilizing the results of computational OT.
Specifically, in \cite{Cuturi2013}, \red{several favorable computational properties of an entropy-regularized version of OT are highlighted. In particular,} entropy-regularized OT problems can be solved efficiently by an iterative algorithm called the Sinkhorn algorithm.
Inspired by this, we propose a dynamical transport algorithm integrating MPC and the Sinkhorn algorithm, which we call {\em Sinkhorn MPC}.
This method incorporates the Sinkhorn iterations into MPC as a dynamic controller and can be seen as simultaneously solving an assignment problem while executing control actions.
The contributions of this paper coming from the introduction of Sinkhorn MPC are as follows:
\begin{enumerate}
	\item[1)] By combining MPC and the Sinkhorn algorithm, the computational effort for determining destinations of agents at each time is reduced substantially;
 	\item[2)] Thanks to the smoothing effect of the entropy regularization, we reveal the global convergence property of Sinkhorn MPC with a \subm{sufficiently} large number of \las{Sinkhorn} iterations;
  \item[3)] For a quadratic control cost, we show the ultimate boundedness and the local asymptotic stability for Sinkhorn MPC without the assumption of the number of \blu{Sinkhorn} iterations.
\end{enumerate}

Compared to a preliminary version of this work~\cite{Ito2022sinkhorn}, the most notable improvement of the current paper is that we consider continuous-time systems and derive 2) the global convergence property for Sinkhorn MPC, which is one of the most crucial properties of dynamical transport algorithms. 
In addition, we provide several illustrative examples, which describe the usefulness of Sinkhorn MPC, and we give the proof of Lemma~\ref{lemma:exponential_convergence}, which is omitted in the preliminary version.

\textit{Organization:}
The remainder of this paper is organized as follows. In Section~\ref{sec:optimal_transport}, we introduce OT between discrete distributions. In Section~\ref{sec:formulation}, we provide the problem formulation.
In Section~\ref{sec:proposed}, we describe the idea of Sinkhorn MPC.
\modi{In Section~\ref{sec:example}, numerical examples illustrate the utility of the proposed method.
Section~\ref{sec:global} is devoted to the global convergence analysis of the proposed method.
In Section~\ref{sec:linear_system}, for a quadratic control cost, we investigate \blu{fundamental} properties of Sinkhorn MPC, \blu{such as local asymptotic stability}.}
Some concluding remarks are given in Section~\ref{sec:conclusion}.

\textit{Notation:}
Let $ \bbR $ denote the set of real numbers.
The set of all positive (resp. nonnegative) vectors in $ \bbR^n $ is denoted by $ \bbR_{> 0}^n $ (resp. $ \bbR_{\ge 0}^n $). We use similar notations for the set of all real matrices $ \bbR^{m\times n} $ and integers $ \bbZ $, respectively.
The set of integers $ \{1,\ldots,N\} $ is denoted by $ \bbra{N} $.
The Euclidean norm is denoted by $ \|\cdot \| $.
For a positive semidefinite matrix $ A $, denote $\|x\|_A := (x^\top A x)^{1/2}$.
The identity matrix of size $n$ is denoted by $I_n$ \subm{or $ I $ when its size is clear in the context}. The matrix norm induced by the Euclidean norm is denoted by $ \| \cdot \|_2 $.
\rev{For vectors $ x_1,\ldots,x_m \in \bbR^n $, a collective vector $ [x_1^\top \ \cdots \ x_m^\top]^\top \in \bbR^{nm} $ is denoted by $ [x_1 ; \ \cdots \ ; x_m] $.}
For $ A = [a_1  \ \cdots \ a_n] \in \bbR^{m\times n}$, we write $ {\rm vec} (A) := \rev{[a_1; \ \cdots \ ;a_n]} $.
For $ \alpha = [\alpha_1 \ \cdots \ \alpha_N ]^\top \in \bbR^N $, the diagonal matrix with diagonal entries $ \{\alpha_i\}_{i=1}^N $ is denoted by \rev{$ \alpha^\diagbox $}.
The block diagonal matrix with diagonal entries $ \{A_i\}_{i=1}^N, A_i \in \bbR^{m\times n} $ is denoted by \rev{$ \{A_i\}_i^\diagbox $}.
Especially when $ A_i = A, \forall i $, $ \{A_i\}_i^\diagbox $ is also denoted by \rev{$  A^{\diagbox,N} $}.
Let $ (\subm{\calM},d) $ be a metric space.
The open ball of radius $ r > 0 $ centered at $ x \in \subm{\calM} $ is denoted by $ B_r(x) := \{y\in \subm{\calM} : d(x,y) < r \} $.
The element-wise division of $ a, b \in \bbR_{>0}^n $ is denoted by $ a \oslash b := [a_1 /b_1 \ \cdots \ a_n/b_n]^\top$.
The $ N $-dimensional vector of ones is denoted by $ \bfone_N $.
The gradient of a function $ f $ with respect to the variable $ x $ is denoted by $ \nabla_x f $.
\rev{For $ x, x' \in \bbR_{> 0}^n $, define an equivalence relation $ \sim $ on $ \bbR_{> 0}^n $ by $ x \sim x' $ if and only if $ \exists r> 0, x = r x' $.}

\section{Background on optimal transport}\label{sec:optimal_transport}
Here, we briefly review OT between discrete distributions $ \mu := \sum_{i=1}^{N} \bfa_i \delta_{x_i}, \nu := \sum_{j=1}^{M} \bfb_j \delta_{y_j} $ where $\bfa \in \Sigma_N := \{ p\in \bbR_{\ge 0}^N : \sum_{i=1}^{N} p_i = 1 \}, \bfb\in \Sigma_M $, $ x_i, y_j \in \bbR^n $, and $ \delta_x $ is the Dirac delta at $ x $.
Given a cost function $ c : \bbR^n \times \bbR^n (\ni (x,y)) \rightarrow \bbR $, which represents the cost of transporting a unit of mass from $ x $ to $ y $, the original formulation of OT due to Monge seeks a map $ \subm{\bfT} : \{x_1,\ldots,x_N\} \rightarrow \{y_1,\ldots,y_M\}$ that solves
\begin{equation}\label{prob:Monge}
	\begin{aligned}
	&\underset{\subm{\bfT}}{\rm minimize} \ \sum_{i \in \bbra{N}} c(x_i, \subm{\bfT}(x_i)) \\
	&\text{subject to} \ \bfb_j = \sum_{i : \subm{\bfT}(x_i) = y_j} \bfa_i, \ \forall j \in \bbra{M} .
	\end{aligned}
\end{equation}
Especially when $ M = N $ and $ \bfa = \bfb = \bfone_N /N $, the optimal map $ \subm{\bfT} $ gives the optimal assignment for transporting agents with the initial states $ \{x_i\}_i $ to the desired states $ \{y_j\}_j $, \red{and then for example, the Hungarian algorithm can be adopted to solve \eqref{prob:Monge}.}
\red{However, this method can be applied only to small problems because it has $ O(N^3) $ complexity.}

On the other hand, the Kantorovich formulation of OT is a linear \subm{program}:
\begin{equation}\label{prob:Kantorovich}
	\underset{P\in \calT(\bfa,\bfb)}{\rm minimize} \ \sum_{i\in \bbra{N},j\in \bbra{M}} C_{ij} P_{ij}
\end{equation}
where $ C_{ij} := c(x_i, y_j) $ and
\[
	\calT(\bfa,\bfb) := \left\{ P \in \bbR_{\ge0}^{N\times M} : P {\bfone}_M = \bfa, \ P^\top \bfone_N = \bfb   \right\}.
\]
A matrix $ P \in \calT(\bfa,\bfb) $, which is called a coupling matrix, represents a transport plan where $ P_{ij} $ describes the amount of mass flowing from $ x_i $ towards $ y_j $.
In particular, when $ M=N $ and $ \bfa = \bfb = {\mathbf 1}_N/N $, there exists an optimal solution from which we can reconstruct an optimal map for \subm{Monge's problem}~\eqref{prob:Monge} \cite[Proposition~2.1]{Peyre2019}.
\red{However, similarly to \eqref{prob:Monge}, for a large number of agents and destinations, the problem~\eqref{prob:Kantorovich} with $ N M $ variables is challenging to solve.}

In view of this, \cite{Cuturi2013} employed entropy regularization to \eqref{prob:Kantorovich}:
\begin{equation}\label{prob:Sinkhorn}
		\underset{P\in \calT(\bfa,\bfb)}{\rm minimize} \ \sum_{i\in \bbra{N},j\in \bbra{M}} C_{ij} P_{ij} - \varepsilon \rmH(P),
\end{equation}
where $ \varepsilon > 0 $ is a regularization parameter and the entropy of $ P $ is defined by $\rmH(P) := - \sum_{i,j} P_{ij} (\log (P_{ij}) - 1)$. 
Define the Gibbs kernel $ K $ associated \subm{with} the cost matrix $ C = (C_{ij}) $ as
\[
	K = (K_{ij}) \in \bbR_{> 0}^{N\times M}, \ K_{ij} := \exp\left(  - C_{ij}/\varepsilon \right) .
\]
Then, a unique solution of the entropic OT problem~\eqref{prob:Sinkhorn} has the form
\begin{equation}\label{eq:opt_form}
	P^* =  (\alpha^*)^\diagbox K  (\beta^*)^\diagbox,
\end{equation}
where the two scaling variables $ (\alpha^*,\beta^*) \in \bbR_{>0}^N \times \bbR_{>0}^M $ are determined by
\begin{equation}\label{eq:opt_scaling_condition}
	\alpha^* = \bfa \oslash [K\beta^*], \ \beta^* = \bfb \oslash [K^\top \alpha^*] .
\end{equation}
The variables $ (\alpha^*, \beta^*) $ can be efficiently computed by the Sinkhorn algorithm:
\begin{equation}\label{eq:Sinkhorn_algorithm}
	\alpha[k+1] = \bfa \oslash [K\beta[k]], \ \beta[k] = \bfb \oslash \left[K^\top \alpha[k]\right]
\end{equation}
where
\[
	\lim_{k\rightarrow \infty} \alpha[k+1]^\diagbox K \beta  [k]^\diagbox  =  P^*, \ \forall \alpha[0] = \alpha_0 \in \bbR_{>0}^N. 
\]

\rev{Now, let us introduce Hilbert's projective metric}
\begin{equation}\label{eq:hilbert_metric}
	\dhil (\beta, \beta') := \log \max_{i,j\in \bbra{M}} \frac{\beta_i \beta'_j}{\beta_j \beta'_i}, \ \beta, \beta' \in \bbR_{>0}^M ,
\end{equation}
\rev{which is a distance on the projective cone $ \bbR_{>0}^M / {\sim} $ (see the Notation in Section~\ref{sec:intro} for $ \sim $)} and is useful for the convergence analysis of the Sinkhorn algorithm; see~\cite[Remark~4.12 \rev{and 4.14}]{Peyre2019}.
Indeed, for any $ (\beta, \beta') \in (\bbR_{> 0}^M )^2 $ and any $ \bar K \in \bbR_{>0}^{N\times M} $, it holds
\begin{equation}\label{ineq:hilbert_metric}
		\dhil (\bar K\beta, \bar K\beta') \le \lambda (\bar K) \dhil (\beta, \beta')
\end{equation}
where
\begin{equation*}
	\lambda (\bar K) := \frac{\sqrt{\eta(\bar K)} -1}{\sqrt{\eta(\bar K)} + 1} < 1,\ \eta(\bar K) := \max_{i,j,k,l} \frac{\bar{K}_{ik} \bar{K}_{jl}}{\bar{K}_{jk} \bar{K}_{il}} .
\end{equation*}
Then it follows from \eqref{ineq:hilbert_metric} that
\begin{align*}
	&\dhil (\beta[k+1], \beta^*) = \dhil \left(\bfb\oslash [K^\top \alpha[k+1]], \bfb\oslash [K^\top \alpha^*]\right) \\
	&=\dhil (K^\top \alpha[k+1], K^\top \alpha^*) \\
	&\le \lambda (K) \dhil (\alpha [k+1], \alpha^*) \le \lambda^2 (K) \dhil (\beta[k], \beta^*) \subm{,}
\end{align*}
which implies $ V_{\subm{\calH}} (\beta) := \dhil (\beta, \beta^*) $ \subm{serves as} a Lyapunov function of \eqref{eq:Sinkhorn_algorithm}, and $ \lim_{k\rightarrow \infty} \beta [k] = \beta^* \in \bbR_{>0}^M / {\sim} $.

\section{Problem formulation}\label{sec:formulation}
In this paper, we consider the problem of \rev{stabilizing agents efficiently to a given discrete distribution over dynamical systems.
This can be formulated as Monge's OT problem.
\begin{problem}
	\label{prob:main_continuous}
	Given initial and desired states $\{x_i^0\}_{i=1}^N, $ $ \{x_j^\sfd\}_{j=1}^N  \in (\mathbb{R}^{n})^N$,
	find control inputs $\{u_i\}_{i=1}^N$ and a permutation $ \sigma : \bbra{N} \rightarrow \bbra{N} $ that solve
	\begin{align}
		&\underset{\sigma }{\rm minimize} ~~ \sum_{i\in \bbra{N}} c_\infty^i (x_i^0, x_{\sigma(i)}^\sfd ) .  \label{eq:infinite_horizon_cost}
	\end{align}
	Here, the cost function $ c_\infty^i : \bbR^n \times \bbR^n \rightarrow \bbR $ is defined by
	\begin{align}
		& &&\hspace{-2cm} c_\infty^i (x_i^0, x_j^\sfd) := \min_{u_i} \ \int_0^\infty  \ell_i (x_i (t), u_i (t) ; x_j^\sfd) \rmd t \label{eq:value_infty}\\
		&\hspace{0.5cm} \text{\rm subject to} &&\red{ \dot{x}_i(t)  = A_i x_i (t)  + B_i u_i (t) }, \label{eq:linear_dynamics_conti}\\
		& &&x_i (t) \in \bbX_i \subseteq \bbR^n , \ \forall t \ge 0 , \label{eq:state_constraint} \\ 
		& &&u_i (t) \in \bbU_i \subseteq \bbR^m , \ \forall t \ge 0 , \label{eq:input_constraint} \\ 
		& &&x_i(0) = x_{i}^0, \label{eq:constraint_initial} \\
		& &&\lim_{t\rightarrow \infty} x_i(t) =  x_j^\sfd,\label{eq:constraint_infty}
	\end{align}
	where $ x_i(t) \in \bbR^n $ denotes the state of the agent $ i $, and $ A_i\in \bbR^{n\times n}, B_i \in \bbR^{n\times m} $.
	\hfill $ \diamondsuit $
\end{problem}
}
Note that the running cost $ \ell_i $ depends not only on the state $ x_i $ and the control input $ u_i $, but also on the destination $ x_{j}^\sfd $. \rev{Throughout this paper, we assume the existence of an optimal solution of OC problems. \red{In addition, we assume that there exists a \blu{constant} input $ \bar{u}_{\blu{ij}} $ under which $ x_i = x_j^\sfd $ is an equilibrium of \eqref{eq:linear_dynamics_conti}.}
\red{A necessary condition for the infinite horizon cost $ c_\infty^i (x_i^0, x_j^\sfd) $ to be finite is that at $ x_i = x_j^\sfd $ \las{and at least one such input} $ u_{i} = \bar{u}_{\blu{ij}} $, there is not a cost incurred, i.e., $ \ell_i(x_j^\sfd, \bar{u}_{\blu{ij}} ; x_j^\sfd) = 0 $.}
\blu{For instance,} if $ B_i $ is square and invertible, $ \bar{u}_{\blu{ij}} =  - B_i^{-1} A_ix_j^\sfd $ makes $ x_i = x_j^\sfd $ an equilibrium.
}

In most cases, the infinite horizon OC problem $ c_\infty^{\magenta{i}} (x_i^0,x_j^\sfd) $ is computationally intractable. To avoid this difficulty, we use MPC, which solves a tractable finite horizon OC problem with a prediction horizon $ T_\rmh > 0 $ at each time:
\begin{align}
	c_{T_\rmh}^{\magenta{i}} (\check{x}_i, x_j^\sfd) := &\min_{u_i} \ \int_{0}^{T_\rmh} \ell_i (x_i(t), u_i (t); x_j^\sfd) \rmd t\label{eq:cost_finite_conti}\\
	&\text{subj. to  \eqref{eq:linear_dynamics_conti}\text{--}\magenta{\eqref{eq:input_constraint}},}~~ x_i(0) = \check{x}_i, \  x_i (T_\rmh) = x_{j}^\sfd . \nonumber
\end{align}
Denote the first control in the optimal sequence of the above problem by $ u_i^{\rm MPC} (\check{x}_i, x_j^\sfd) $.
Also for $ \check{x} = [\check{x}_1 ; \cdots ; \check{x}_N] \in \bbR^{nN} $, denote by $ \sigma(\cdot; \check{x}) $ the optimal permutation of the following problem:
\begin{align}
	&\underset{\sigma }{\rm minimize} ~~ \sum_{i\in \bbra{N}} c_{T_\rmh}^i (\check{x}_i, x_{\sigma(i)}^\sfd ) .  \label{eq:finite_horizon_OT}
\end{align}
Then the dynamics \eqref{eq:linear_dynamics_conti} under MPC for Problem~\ref{prob:main_continuous} is given by
\begin{equation}\label{eq:permutation_dynamics}
  \dot{x}_i (t) = A_i x_i (t) + B_i u_i^{\rm MPC} \left(x_i(t), x_{\sigma(i;x(t))}^\sfd \right) , \ \forall i\in \bbra{N} ,
\end{equation}
where $ x(t) := [x_1(t);\cdots ; x_N(t)] $.
Note that along the trajectory $ x(t) $, at several times, the permutation $ \sigma(\cdot;x(t)) $ changes in general.
The state-dependent permutation $ \sigma(\cdot; x(t)) $ is expected to \blu{reduce the cost accumulated during the transport} \las{more} than the permutation $ \sigma(\cdot; x^0) $ that is fixed at the initial time $ t=0 $.
Despite the merit, the state-dependency of the permutation poses the following computational and theoretical difficulties:
\begin{itemize}
	\item Solving the assignment problem~\eqref{eq:finite_horizon_OT} at each time leads to the high computational burden when $ N $ is large;
 	\item The optimal permutation $ \sigma(\cdot; x) $ is not continuous in $ x $. That is, the target states $ \{x_{\sigma(i;x(t))}^\sfd \}_i $ for the agents change discontinuously along the trajectory $ x(t) $, and this makes \blu{it} difficult to ensure the \blu{convergence} of the dynamics~\eqref{eq:permutation_dynamics}.
\end{itemize}
In the remainder of this paper, we reveal that entropy regularization mitigates the above problems.

\section{MPC with entropy-regularized optimal transport}\label{sec:proposed}
\subsection{Introduction of the entropy regularization to MPC}
\begin{figure}[tb]
	\centering
	\includegraphics[scale=0.4]{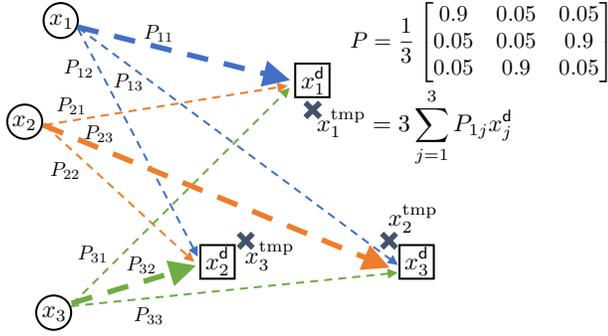}
	\caption{\comm{\subm{Three} agents $ \{x_i\} $, desired states $ \{x_{\subm{j}}^\sfd\} $, and temporary targets $ \{x_i^\tmp (P)\} $ determined by the barycentric projection~\eqref{eq:barycentric_target} for a given coupling matrix $ P $.}}
	\label{fig:image}
\end{figure}

Now, to avoid the issues observed in the previous section, we employ the entropy regularization.
To this end, we first consider the linear program:
\begin{equation}\label{eq:equi_LP}
	\underset{P\in \calT(\comm{\frac{\one}{N},\frac{\one}{N}})}{\rm minimize} \ \sum_{i,j\in \bbra{N}} C_{ij} (x) P_{ij},
\end{equation}
where $ C_{ij} (x) := c_{T_\rmh}^i (x_i,x_j^\sfd), \ x = [x_1;\cdots;x_N] $. Then as mentioned in Section~\ref{sec:optimal_transport}, the optimal permutation $ \sigma $ can be obtained by the optimal permutation matrix $ P^\sigma $ of \eqref{eq:equi_LP} satisfying $ P_{ij}^\sigma = 1/N $ if $ j= \sigma(i) $, and $ 0 $, otherwise.
Next, we introduce the entropy regularization to \eqref{eq:equi_LP} as in \eqref{prob:Sinkhorn}.
Then, based on the optimal coupling $ P^* $ of the entropic OT problem, we determine a target state for each agent.
Specifically, we introduce a map $ x_i^\tmp : \bbR_{\ge 0}^{N\times N} \rightarrow \calX (\subset \bbR^n) $ as a policy to determine a temporary target $ x_i^\tmp (P^*) $ for agent $ i $.
\comm{We call $ x_i^\tmp $ a navigator function.}
A typical \comm{navigator function} to approximate Monge's OT map from a coupling matrix $ P $ is the so-called barycentric projection~\cite[Remark~4.11]{Peyre2019}:
\begin{equation}\label{eq:barycentric_target}
x_i^\tmp (P) = N \sum_{j=1}^{N} P_{ij} x_j^\sfd, \ P \in \bbR_{\ge 0}^{N\times N} .
\end{equation}
Note that, for a permutation matrix $ P^\sigma $, it holds $ N\sum_{j=1}^N P_{ij}^\sigma x_j^\sfd = x_{\sigma(i)}^\sfd $.
\comm{Fig.~\ref{fig:image} illustrates the states of three agents $ \{x_i\} $, destinations $ \{x_{\subm{j}}^\sfd\} $, and \blu{temporary} targets $ \{x_i^\tmp(P)\} $ determined by the barycentric projection~\eqref{eq:barycentric_target} for a given coupling matrix $ P $.}

Now, we propose to use the control law \las{$ u_i^\mpc \bigl(x_i(t), $ $ x_i^\tmp (P^*(x(t)))\bigr) $} where
\begin{equation}
	P^* (x) := \argmin_{P\in \calT(\frac{\one}{N}, \frac{\one}{N})} \ \sum_{i,j\in \bbra{N}} C_{ij} (x) P_{ij} - \varepsilon \rmH(P). \label{eq:P_opt_sink}
\end{equation}
In summary, for any given \subm{navigator function}~$ x_i^\tmp $ and $ \varepsilon > 0 $, the dynamics of the agents are written as
\begin{align}
  &\dot{x}_i (t) = A_i x_i (t) + B_i u_i^{\rm MPC} \left(x_i(t), x^{\rm tmp}_{i} \Bigl(P^*(x(t))\Bigr) \right) , \nonumber \\
  &\hspace{6cm} \forall i\in \bbra{N} , \label{eq:mpc_dynamics_conti}\\
  &x_i (0) = x_i^0, \ \forall i \in \bbra{N} . \nonumber 
\end{align}
The entropy regularization enables to use the Sinkhorn algorithm~\eqref{eq:Sinkhorn_algorithm}, which contributes to \blu{reducing} the computational burden of determining target states at each time.
In addition, we will see that the entropy regularization also enables to analyze the global convergence property of \eqref{eq:mpc_dynamics_conti} in Section~\ref{sec:global}.

\subsection{Integrating MPC and the Sinkhorn algorithm}\label{subsec:proposed}
In the previous subsection, \subm{it was implicitly} assumed that at each time, the optimal coupling $ P^*(x(t)) $ is available for determining temporary targets for agents.
The Sinkhorn algorithm achieves a speed-up in the computation of an optimal coupling.
However, in challenging situations in which the number of agents is very large and the sampling time is small, only a few Sinkhorn iterations are allowed. In such a case, an available approximate coupling matrix may not be close enough to the optimal coupling, and therefore the above assumption is not valid.
To address this issue, in this subsection, we propose to integrate MPC and the Sinkhorn algorithm.
Since the Sinkhorn algorithm works in discrete time, we consider a time-discretized version of \eqref{eq:linear_dynamics_conti} obtained by e.g., a zero-order hold discretization\footnote{Throughout this paper, we use bold symbols for discrete-time systems while we use italic letters for continuous-time systems.}:
\begin{equation}\label{eq:linear_dynamics_disc}
	\bmx_i [k+1]  = \bmA_i \bmx_i [k]  + \bmB_i \bmu_i [k], \ k\in \bbZ_{\ge 0} .
\end{equation}
Then, the cost function $ \bmc_{\tau_\rmh}^i $ with a finite horizon $ \tau_\rmh \in \bbZ_{>0} $ is defined by
\begin{align}
	& &&\hspace{-1cm} \bmc_{\tau_\rmh}^i (\check{x}_i, x_j^\sfd) := \min_{\bmu_i} \ \sum_{k=0}^{\tau_\rmh-1}  {\boldsymbol \ell}_i \left(\bmx_i [k], \bmu_i [k] ; x_j^\sfd \right) \label{eq:cost_disc}\\
	&\hspace{-0.1cm} \text{\subm{subj. to}} &&\eqref{eq:linear_dynamics_disc},\nonumber\\
	& &&\bmx_i[k] \in \bbX_i, \ \bmu_i [k] \in \bbU_i , \ \forall k \in \subm{\bbra{\tau_\rmh -1}\cup \{0\}} , \nonumber \\ 
	& &&\bmx_i[0] = \check{x}_{i}, \label{eq:constraint_initial_disc} \\
	& &&\bmx_i[\tau_\rmh] =  x_j^\sfd. \label{eq:constraint_end_disc}
\end{align}
Denote the first control in the optimal sequence of the above problem by $ \bmu_i^{\rm MPC} (\check{x}_i, x_j^\sfd) $. Let $ \bmx[k] := [\bmx_1[k];\cdots; \bmx_N[k]] $ and
\begin{equation}
	\bmP^* (x) := \argmin_{P\in \calT(\frac{\one}{N}, \frac{\one}{N})} \ \sum_{i,j\in \bbra{N}} \bmC_{ij} (x) P_{ij} - \varepsilon \rmH(P), \label{eq:P_opt_sink_disc}
\end{equation}
where $ \bmC_{ij} (x) := \bmc_{\tau_\rmh}^i (x_i, x_j^\sfd) $.
Note that if $ \bmc_{\tau_\rmh}^i (x_{\blu{i}},x_j^\sfd) $ is continuous in $ x_{\blu{i}} $ for all $ i\in \bbra{N} $, then from the relations~\eqref{eq:opt_form},~\eqref{eq:opt_scaling_condition}, $ \bmP^* $ is continuous.
Hence, it is expected that if we take a sampling time for \eqref{eq:linear_dynamics_disc} appropriately so that the difference between $ \bmx_i [k+1] $ and $ \bmx_i [k] $ is small, then the scaling variables $ (\alpha^*[k+1], \beta^*[k+1]) $ for $ \bmP^*(\bmx[k+1]) $ are close to the variables $ (\alpha^*[k],\beta^*[k]) $ for $ \bmP^* (\bmx[k]) $. This implies that $ (\alpha^*[k],\beta^*[k]) $ yield good initial estimates of $ (\alpha^*[k+1], \beta^*[k+1]) $.

Based on this observation, we present a dynamical transport algorithm integrating MPC and the Sinkhorn algorithm.
Let $ S[k] \in \bbZ_{>0} $ be the number of \blu{Sinkhorn} iterations at time $ k $.
\rev{For any given \comm{navigator function} $ x_i^\tmp $ and $ \varepsilon > 0 $,} the proposed algorithm, which we call {\em Sinkhorn MPC} is given as the following dynamics where the Sinkhorn algorithm behaves as a dynamic controller.\\
{\bf Sinkhorn MPC:}
\begin{align}
&\bmx_i[k+1] = \bmA_i \bmx_i [k] + \bmB_i \bmu_i^{\rm MPC} \bigl(\bmx_i [k], x_i^\tmp \left(P[k] \right)  \bigr),\nonumber\\
&\hspace{6cm} \forall i \in \bbra{N}, \label{eq:SMPC_x} \\
&P[k] = \alpha\left[k,S[k]+1\right]^\diagbox \bmK(\bmx[k]) \beta\left[k,{S}[k]\right]^\diagbox, \label{eq:SMPC_P} \\
&\text{Sinkhorn iterations:} \nonumber\\
&\begin{cases}
\alpha\left[k, {l+1}\right] = \bfone_N/N \oslash \left[\bmK(\bmx[k])\beta[k,l] \right], \\
\beta[k, l] = \bfone_N/N \oslash \left[\bmK(\bmx[k])^\top \alpha[k,l] \right] , 
\end{cases} \ l \in \bbra{S[k]}, \label{eq:SMPC_sink}\\
&\alpha[k+1,1] = \alpha[k,S[k]+1],  \label{eq:SMPC_next} \\ 
&\bmx_i [0] = x_i^0, \ \alpha[0,1] = \alpha_0, \nonumber 
\end{align}
where
\begin{align*}
\bmK_{ij}(x) := \exp\left(- \frac{\bmc_{\tau_\rmh}^{\magenta{i}} (x_i, x_{j}^\sfd )}{\varepsilon} \right), \ x = [x_1;  \cdots  ;x_N] ,
\end{align*}
and the initial value \las{$ \alpha_0 \in \bbR_{>0}^N $} is arbitrary.
\hfill $ \diamondsuit $

The important point here is that for the initial value ${ \alpha[k+1,1] }$ of the Sinkhorn iterations at time $ k+1 $, we use the final value $ \alpha[k,{S[k] + 1}] $ of the iterations at the previous time $ k $.
When the difference between $ \bmx[k+1] $ and $ \bmx[k] $ is small, $ \alpha[k,{S[k] + 1}] $ will be a good initial estimate of the scaling variable $ \alpha^*[k+1] $ for $ \bmP^*(\bmx[k+1]) $ even if $ S[k] $ is small.
A convenient way to determine $ S[\cdot] $ is to fix it to a suitable constant $ \hat{S}\in \bbZ_{>0} $ in terms of computation time.

The pseudocode of Sinkhorn MPC is described in Algorithm~\ref{algo:sinkhorn_mpc}.
Note that, of course, the proposed method can be applied to systems that are not discretizations of continuous-time systems and originally evolve in discrete time.

\begin{algorithm}
	\caption{Sinkhorn MPC}\label{algo:sinkhorn_mpc}
	\begin{algorithmic}[1]
		\renewcommand{\algorithmicrequire}{\textbf{Input:}}
		\renewcommand{\algorithmicensure}{\textbf{Output:}}
		\Require $ \{x_i^0\}_{i=1}^N, \{x_j^\sfd\}_{j=1}^N, \alpha_0, \ \varepsilon, \ \tau_\rmh, \ S[\cdot] $\\
		\textit{Initialization}: $ \bmx_i[0] := x_i^0,  \forall i, \ \alpha := \alpha_0$
		\For{$ k = 0, 1, 2,\ldots$}
		\For{$ i,j \in \bbra{N}$ ({\it run in parallel})}
		\State Compute $ \bmc_{\tau_\rmh}^i (\bmx_i[k], x_j^\sfd) $ defined in \eqref{eq:cost_disc}
		\State $ \bmK_{ij} := \exp\left(- \frac{\bmc_{\tau_\rmh}^{\magenta{i}} (\bmx_i [k], x_{j}^\sfd )}{\varepsilon} \right) $
			\EndFor
			\For{$ l = 1,2,\ldots,S[k] $}
			\State $ \beta := \bfone_N/N \oslash [\bmK^\top \alpha] $
			\State $ \alpha := \bfone_N/N \oslash [\bmK \beta] $
		\EndFor
		\State $ P: = \alpha^\diagbox \bmK \beta^\diagbox $
		\For{$ i \in \bbra{N} $ ({\it run in parallel})}
			\State Compute $ \bmu_i^\mpc (\bmx_i[k], x_i^\tmp (P)) $ and apply it to the agent $ i $
		\EndFor
		\EndFor
	\end{algorithmic}
\end{algorithm}

\section{Illustrative examples}\label{sec:example}
\subsection{Quadratic cost}\label{subsec:quad_invertible}
This section gives examples for Sinkhorn MPC.
First we consider \blu{a} quadratic cost
\begin{equation}\label{eq:energy}
	\bmell_i (x_i, u_i; x_j^\sfd) = \|u_i - \bmB_i^{-1} (x_j^\sfd - \bmA_i x_j^\sfd)\|^2 ,
\end{equation}
where we assumed the invertibility of $ \bmB_i $, and let $ {\bbX_i = \bbU_i = \bbR^n }$.
Note that for a constant input $ \bar{\bmu}_{\blu{ij}} := \bmB_i^{-1} (x_j^\sfd - \bmA_i x_j^\sfd) $, under which $ \bmx_i = \bmx_j^\sfd $ is an equilibrium of \eqref{eq:linear_dynamics_disc}, it holds $ \bmell_i (x_j^\sfd, \bar{\bmu}_{\blu{ij}} ; x_j^\sfd) = 0 $.
Then, the dynamics under Sinkhorn MPC can be written as follows \cite[\rev{Section~2.2, pp.~37-39}]{Lewis2012}:
\begin{align}
&\bmx_i[k+1] = \bar{\bmA}_i \bmx_i [k] +  (I - \bar{\bmA}_i) x_i^\tmp (P[k]), \label{eq:SMPC_linear_x} \\
&\bmu_i^{\rm MPC} (x_i, \red{\hat{x}}) = - \bmB_i^\top (\bmA_i^\top)^{\tau_\rmh - 1} \bmG_{i,\tau_\rmh}^{-1} \bmA_i^{\tau_\rmh} (x_i - \red{\hat{x}} ) \nonumber\\
&\hspace{2.3cm} + \bmB_i^{-1} ( \red{\hat{x}} - \bmA_i \red{\hat{x}}), \ \forall i \in \bbra{N}, \ \red{\forall \subm{x_i}, \hat{x}\in \bbR^n} \nonumber
\end{align}
with \eqref{eq:SMPC_P},~\eqref{eq:SMPC_sink} where
\begin{align*}
&\bmK_{ij}(x)  = \exp \left(- \frac{  \|x_i - x_j^\sfd\|_{\calG_{i}}^2 }{\varepsilon}   \right), \\
&\calG_{i} := (\bmA_i^{\tau_\rmh})^\top \bmG_{i,\tau_\rmh}^{-1}  \bmA_i^{\tau_\rmh}, \ \bmG_{i,\tau_\rmh} := \sum_{k=0}^{\tau_\rmh -1} \bmA_i^k \bmB_i \bmB_i^\top (\bmA_i^\top)^k, \\
&\bar{\bmA}_i := \bmA_i - \bmB_i \bmB_i^\top (\bmA_i^\top)^{\tau_\rmh -1} \bmG_{i,\tau_\rmh}^{-1} \bmA_i^{\tau_\rmh} .
\end{align*}
In the examples below, we use the barycentric target~\eqref{eq:barycentric_target} \subm{as a navigator function}.

\begin{figure}[tb]
	\centering
	\includegraphics[scale=0.45]{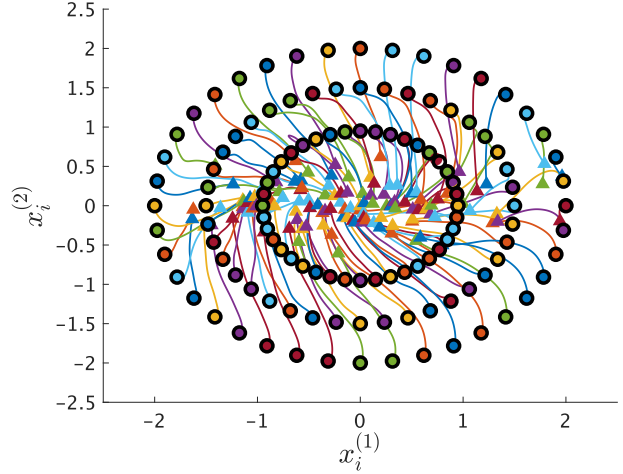}
	\caption{Trajectories $ \bmx_i[k] = [\bmx_i^{(1)}[k] \ \bmx_i^{(2)} [k]]^\top $ of 120 agents for \eqref{eq:ex2d} (solid), initial states (filled triangles), steady states (filled circles), and desired states (black circles).}
	\label{fig:2d_trajectory}
\end{figure}

\begin{figure}[tb]
	\centering
	\includegraphics[scale=0.35]{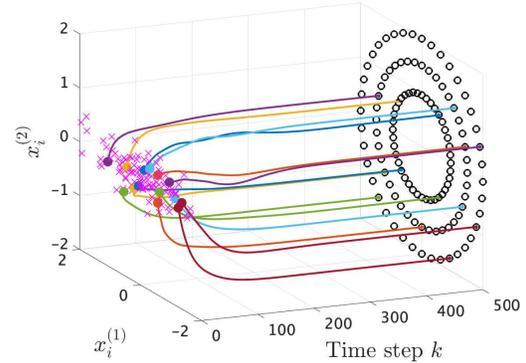}
	\caption{Trajectories of $ 14 $ agents out of $ 120 $ agents for \eqref{eq:ex2d} (solid), initial states (magenta crosses), and desired states (black circles).}
	\label{fig:3d_trajectory}
\end{figure}

\begin{figure}[tb]
	\centering
	\includegraphics[scale=0.45]{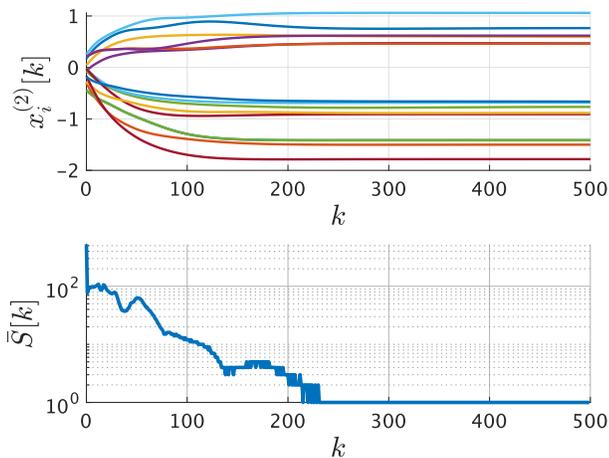}
	\caption{Trajectories $ \bmx_i^{(2)} [k] $ of 14 agents (top) and the number of \blu{Sinkhorn} iterations at each time $ \bar{S}[k] $ (bottom, \res{semi-log plot}) for \eqref{eq:ex2d}.}
	\label{fig:iteration}
\end{figure}

First, consider \eqref{eq:linear_dynamics_conti} with
\begin{equation}\label{eq:ex2d_conti}
	A_i = 
	\begin{bmatrix}
		2 & 1.3\\
		-0.5 & 1
	\end{bmatrix}, \
	B_i = I_2, \ \forall i\in \bbra{N} .
\end{equation}
By using the Euler method with a step size $ 0.02 $, we obtain 
\begin{equation}\label{eq:ex2d}
	\bmA_i = 
	\begin{bmatrix}
		1.04 & 0.026\\
		-0.01 & 1.02
	\end{bmatrix}, \
	\bmB_i = 0.02 I_2, \ \forall i\in \bbra{N} .
\end{equation}
Set $ N = 120, \ \varepsilon = 2.0, \ \tau_\rmh = 100, \ \alpha_0 = \bfone_N $.
Here, we consider the case where the optimal coupling $ \bmP^* $ is available at each time.
Specifically, rather than using $ S[\cdot] $ fixed beforehand for Sinkhorn MPC, we employ the stopping criterion for the Sinkhorn iterations~\cite[Remark~4.14]{Peyre2019} given by
\begin{equation}\label{ineq:stop}
	\| P[k] \bfone_N - \bfone_N/N \|_1 + \|P[k]^\top \bfone_N - \bfone_N/N \|_1 < 0.005,
\end{equation}
where $ \| \cdot \|_1 $ denotes the $ \ell^1 $-norm.
For given initial and desired states, the trajectories of the agents governed by \eqref{eq:SMPC_linear_x} with \eqref{eq:SMPC_P}--\eqref{eq:SMPC_next},~\subm{\eqref{eq:ex2d}}, \eqref{ineq:stop} are illustrated in Figs.~\ref{fig:2d_trajectory},~\ref{fig:3d_trajectory}.
It can be seen that the agents converge sufficiently close to the target states. We will study the convergence property in Sections~\ref{sec:global},~\ref{sec:linear_system}.
The number of \blu{Sinkhorn} iterations $ \bar{S}[k] $ satisfying \eqref{ineq:stop} at each time $ k $ is shown in Fig.~\ref{fig:iteration}. The number of \blu{iterations} is drastically reduced from $ \bar{S}[0] \simeq 520 $ to $ \bar{S}[1] \simeq 100 $ in one time step, and $ \bar{S} $ continues to decrease as $ k $ increases.
This clarifies that the optimal scaling variables $ (\alpha^*[k],\beta^*[k]) $ can be used for good initial estimates of $ (\alpha^*[k+1], \beta^*[k+1]) $ as expected in Subsection~\ref{subsec:proposed}.

The computation time for one Sinkhorn iteration and the number of \blu{Sinkhorn} iterations~$ \bar{S}[0] $ at the initial time with different $ N $ are shown in Table~\ref{table:computation}. The algorithm has been implemented in MATLAB on MacBook~Pro with Apple~M1~Pro.
Table~\ref{table:computation} also shows the computation time for solving an optimal assignment problem to obtain \subm{the permutation} $ \sigma(\cdot;\bmx[k]) $ by the Hungarian algorithm~\cite{Cao2022}.
As can be seen, the Hungarian algorithm is not scalable and thus not suitable for MPC.
Hence, introducing the entropy regularization to MPC contributes to reducing the computational burden.

\renewcommand{\arraystretch}{1.3}
\begin{table*}[hbtp]
  \caption{Computation time for one Sinkhorn iteration and the Hungarian algorithm, and the number of \blu{iterations} $ \bar{S}[0] $.}
  \label{table:computation}
  \centering
  \begin{tabular}{lrrrr}
    \hline
      & \begin{tabular}{c}Computation time for\\one Sinkhorn iteration \end{tabular}  & \begin{tabular}{c}$ \bar{S}[0] $ \\with $ \varepsilon = 2.0 $ \end{tabular}  &\begin{tabular}{c} $ \bar{S}[0] $\\with $ \varepsilon = 4.0 $ \end{tabular} & Hungarian algorithm~\cite{Cao2022}\\
    \hline 
    $ N  =120 $  & $ 0.0060 $~ms  & $ 527 $  & $ 267 $ & $ 0.1 $~s\\
    $ N = 500 $  & $ 0.21 $~ms   & $ 381 $ & $ 189 $ & $ 0.6 $~s\\
		$ N = 1000 $  & $ 0.83 $~ms   & $ 379 $ & $ 186 $ & $ 3.0 $~s\\
		$ N = 3000 $  & $ 3.0 $~ms   & $ 256 $ &  $ 129 $&  $ 57.0 $~s\\
    \hline
  \end{tabular}
\end{table*}
\renewcommand{\arraystretch}{1.0}

\subsection{Effect of the number of \blu{Sinkhorn} iterations and the regularization parameter on Sinkhorn MPC}
\begin{figure}[tb]
	\centering
	\includegraphics[scale=0.35]{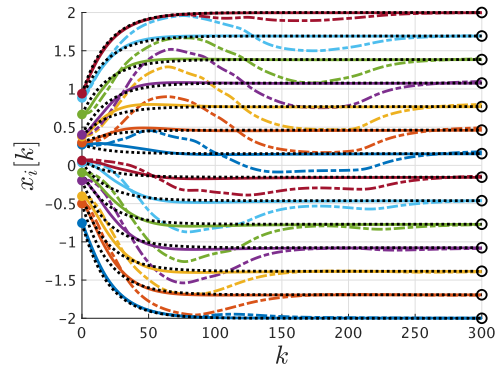}
	\caption{Trajectories $ \bmx_i [k] $ of 14 agents for \eqref{eq:ex1d} with $ S[k] \equiv 1 $ (chain) and $ S[k] \equiv 5 $ (solid), respectively, and desired states (black circles). The black dotted lines indicate trajectories without \blu{the} regularization ($ \varepsilon = 0 $) following \eqref{eq:permutation_dynamics}.}
	\label{fig:iteration_compare}
\end{figure}

Here, we describe how the number of \blu{Sinkhorn} iterations $ S $ affects the bahavior of Sinkhorn MPC.
To this end, consider a simple case for~\eqref{eq:energy} with $ N = 14, \ \tau_\rmh = 20, \ \varepsilon = 0.1 $, and
\begin{equation}\label{eq:ex1d}
	\bmA_i = 1, \ \bmB_i = 0.1, \ \forall i \in \bbra{N}.
\end{equation}
Then the trajectories of the agents with $ S[k] = 1, \forall k $ and $ S[k] = 5, \forall k $ are illustrated in Fig.~\ref{fig:iteration_compare}. 
Also, the trajectories without the regularization $ (\varepsilon = 0) $ following the discretized version of \eqref{eq:permutation_dynamics} are shown with the black dotted lines.
In this example, for all the cases, $ \bmx[k] $ converges to almost the same point close enough to the desired distribution.
The total energy cost $ \sum_{i,k} \|\bmu_i[k] \|^2 $ for the case without the regularization is $ 21.6 $.
It can be seen from Fig.~\ref{fig:iteration_compare} that one iteration per time step is not enough to determine an appropriate destination for each agent while performing control and results in the total energy cost $ 34.7 $.
On the other hand, the trajectories for five iterations are similar to the trajectories without the regularization, and the total energy cost is $ 19.1 $. Note that since we use MPC, the total cost for Sinkhorn MPC can be smaller than for the case without the regularization as in this example.

\begin{figure}[t]
	\begin{minipage}[b]{0.53\linewidth}
		\hspace{-0.45cm}
		\includegraphics[keepaspectratio, scale=0.3]
		{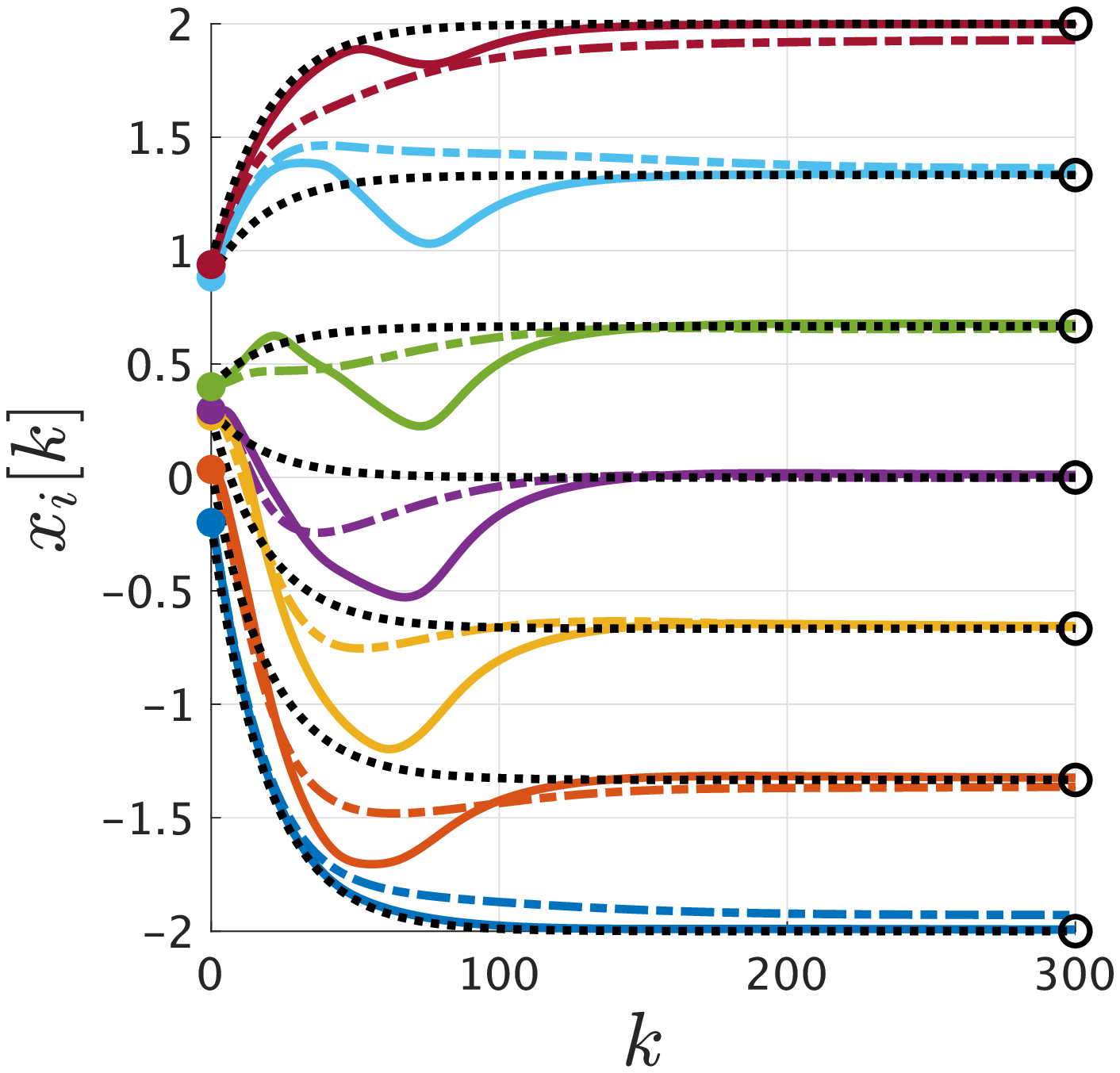}
		\subcaption{}\label{fig:ex1d}
	\end{minipage}
	\begin{minipage}[b]{0.46\linewidth}
		\includegraphics[keepaspectratio, scale=0.3]
		{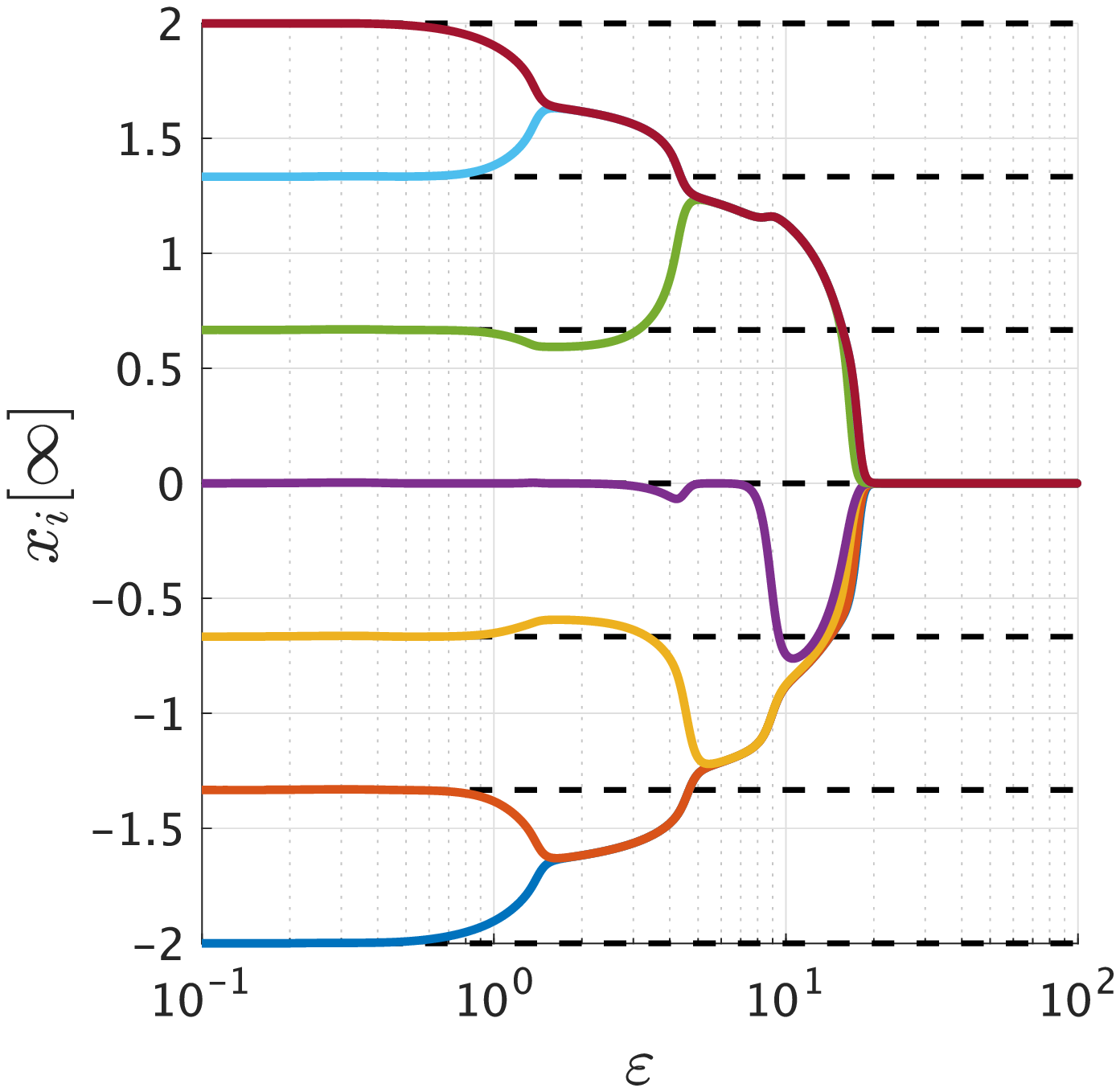}
		\subcaption{}\label{fig:steady_state}
	\end{minipage}

	\caption{(a) Trajectories $ \bmx_i[k] $ of 7 agents for \eqref{eq:ex1d} with $ \varepsilon = 0.9 $ (chain) and $ \varepsilon = 0.4 $ (solid), respectively, \rev{trajectories without the regularization (dotted)}, and desired states (black circles).  \\
	\rev{\subm{(b)} Semi-\blu{log} plot of steady states $ \lim_{k\rightarrow\infty} \bmx_i [k] $ with different $ \varepsilon \in [10^{-1},10^2] $ (solid) and desired states (dashed).}}
\end{figure}

Next, we investigate the effect of the regularization parameter $ \varepsilon $ on Sinkhorn MPC.
We continue to consider \eqref{eq:ex1d} and set $ N = 7, S[k] = 1,  \forall k $.
Then the trajectories of the agents with $ \varepsilon = 0.4, 0.9 $ are shown in Fig.~\ref{fig:ex1d}.
As can be seen, the overshoot/undershoot is reduced for larger $ \varepsilon $ while the limiting values of the states deviate from the desired states.
In other words, the parameter $ \varepsilon $ reflects the trade-off between the stationary and transient behaviors of the dynamics under Sinkhorn MPC.
In fact, it is known that the convergence of the Sinkhorn algorithm deteriorates as $ \varepsilon \rightarrow + 0 $~\cite[Remark~4.14]{Peyre2019}; see Table~\ref{table:computation}. This degrades the transient behaviors under Sinkhorn MPC.
Taking larger $ S[k] $ remedies this issue.
The steady states $ \lim_{k\rightarrow \infty} \bmx_i [k] $ for a fixed set of initial states as a function of $ \varepsilon $ are illustrated in Fig.~\ref{fig:steady_state}. Note that for different initial states, we obtained the same result or the one flipped upside down.
The obtained behavior is due to the fact that as $ \varepsilon $ becomes larger, the optimal coupling of the entropic OT problem is more blurred to the maximum entropy coupling $ \bfone_N \bfone_N^\top /N^2 $~\cite[Proposition~4.1]{Peyre2019}.
The behavior of the equilibrium points for Sinkhorn MPC as $ \varepsilon \rightarrow + 0 $ will be revealed in Lemma~\ref{lemma:exponential_convergence} in Section~\ref{sec:linear_system}.
Although we have considered the simple setting~\eqref{eq:ex1d}, the above observations apply to the general case.

\subsection{Non-quadratic cost}
\begin{figure}[tb]
	\centering
	\includegraphics[scale=0.35]{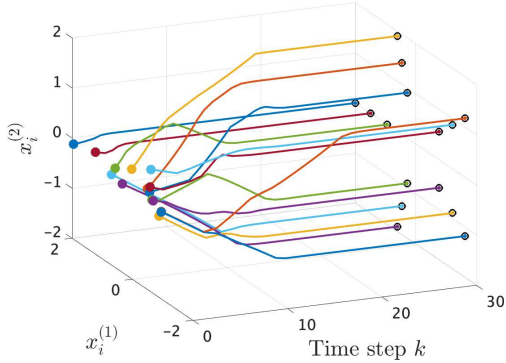}
	\caption{Trajectories of $ 15 $ agents for \eqref{eq:l1_cost} with $ \varepsilon = 2.0 $ (solid) and desired states (black circles).}
	\label{fig:lasso}
\end{figure}

Lastly, we investigate the behavior of Sinkhorn MPC for \las{the} non-quadratic cost and \las{the} bounded input spaces:
\begin{align}
	&\bmell_i (x_i, u_i; x_j^\sfd) = \|x_i - x_j^\sfd \|^2  +  \|u_i - \bmB_i^{-1} (x_j^\sfd - \bmA_i x_j^\sfd)\|_1 \subm{,}  \label{eq:l1_cost} \\
	&\bbU_i =  \subm{\left\{ u_i = [u_i^{(1)} \ u_i^{(2)}]^\top \in \bbR^{2} : |u_i^{(j)}| \le 5,\  j = 1,2 \right\} } , \nonumber\\
	&\subm{\bbX_i = \bbR^2}, \hspace{5cm} \forall i\in \bbra{N} .\nonumber
\end{align}
Here $ \ell^1 $-norm encourages $ u_i $ to be identically equal to $ \bmB_i^{-1} (x_j^\sfd - \bmA_i x_j^\sfd) $ for long time and is used for sparse optimal control~\cite{Nagahara2016,Ito2021sparse}.
For the computation of $ \bmc_{\tau_\rmh} $ and $ \bmu_i^\mpc $, we used {\tt cvx} package with MATLAB~\cite{cvx}.
Set $ N  = 15, \varepsilon = 2.0, \tau_\rmh = 50, S [k] = 20, \forall k $ and consider \eqref{eq:ex2d}.
Then, Fig.~\ref{fig:lasso} shows the trajectories under Sinkhorn MPC with the barycentric target~\eqref{eq:barycentric_target}. As can be seen, the agents have achieved the desired distribution.
This result shows that Sinkhorn MPC achieves the transport also for the non-quadratic cost.

\section{Global convergence property of Sinkhorn MPC}\label{sec:global}
In the remainder of this paper, we investigate the fundamental properties of Sinkhorn MPC.
In this section, we consider the case where the sampling time for obtaining \eqref{eq:linear_dynamics_disc} is small and $ S $ is large enough so that the dynamics~\eqref{eq:mpc_dynamics_conti} is well approximated by the discretized system~\eqref{eq:SMPC_x}.
Then, we analyze the global convergence property of the dynamics~\eqref{eq:mpc_dynamics_conti}.
To this end, we suppose the following condition holds.
\begin{assumption}\label{ass:positive_cost}
  For all $ i\in \bbra{N} $, $ \bbX_i = \bbU_i = \bbR^n $, $ B_i $ is invertible, and
	the function $ \ell_i $ satisfies
  \begin{equation}\label{eq:positive_cost}
    \ell_i (x_i, u_i; x_j)
    \begin{cases}
      =0  & \text{if} \ x_i =x_j \ \text{and} \ u_i = -B_i^{-1} A_i x_i, \\
      >0  & \text{otherwise}.
    \end{cases}
  \end{equation}
  \hfill $ \diamondsuit $
\end{assumption}
Then, a point $ x^\sfe = [x_1^\sfe;\cdots;x_N^\sfe] \in \bbR^{nN} $ satisfying
\begin{equation}\label{eq:equi_general}
  x_i^\sfe = x^{\rm tmp}_{i} (P^*(x^\sfe)), \ \forall i \in \bbra{N}
\end{equation}
is an equilibrium of \eqref{eq:mpc_dynamics_conti}.
Indeed, by \eqref{eq:positive_cost} and \eqref{eq:equi_general}, the constant input $ u_i (t) \equiv -B_i^{-1} A_i x_i^\sfe $, under which the state process starting from $ x_i(0) = x_i^\sfe $ is the constant $ x_i(t) \equiv x_i^\sfe $, is the unique optimal solution of the OC problem \eqref{eq:cost_finite_conti} with $ \check{x}_i = x_i^\sfe, x_j^\sfd  = x^{\rm tmp}_{i} (P^*(x^\sfe)) = x_i^\sfe $. 
Therefore,
\[
  u_{i}^\mpc \left(x_i^\sfe, x_i^\tmp (P^*(x^\sfe)) \right) = -B_i^{-1} A_i x_i^\sfe ,
\]
and
\[
  A_i x_i^\sfe + B_i u_i^\mpc \left(x_i^\sfe, x_i^\tmp (P^*(x^\sfe))\right) = 0 .
\]

The following proposition ensures the existence of a point satisfying \eqref{eq:equi_general}.
\begin{proposition}\label{prop:fixed_point_general}
  Assume that the codomain $ \calX $ of $ x_i^\tmp $ is a compact convex set and for all $ i\in \bbra{N} $, $ x_i^\tmp : \bbR_{\ge 0}^{N\times N} \rightarrow \calX $ is continuous.
	Assume further that for all $ i,j\in \bbra{N} $, $ c_{T_\rmh}^i (x_i,x_j^\sfd) $ is continuous in $ x_i $.
  Then, the set ${ \calR := \{ x^\sfe \in \bbR^{nN} : x_i^\sfe =  x^{\rm tmp}_{i} (P^*(x^\sfe)), \ \forall i \in \bbra{N}\} }$ is non-empty.
\end{proposition}
\begin{pf}
	Define a map $ h : \bbR^{nN} \rightarrow \bbR^{nN} $ as
	\begin{equation}
	h(x) := \left[x_1^\tmp (P^*(x)); \cdots ; x_N^\tmp (P^*(x))\right]
	, \ x \in \bbR^{nN} . \label{eq:func_fixed_point}
	\end{equation}
	It is obvious from \eqref{eq:opt_form},~\eqref{eq:opt_scaling_condition}, and the continuity of $ c_{T_\rmh}^i $ that $ P^* $ is continuous.
  Since $ x_i^\tmp $ and $ P^* $ are continuous, $ h $ is also continuous.
	The set of all fixed points of $ h $ coincides with $ \calR $.
	For brevity, we abuse notation and regard $ \calX^N $ as a subset of $ \bbR^{nN} $.
	Let $ h_{\calX} : {\calX}^N \rightarrow {\calX}^N $ be the restriction of $ h $ to $ {\calX}^N $.
	Now we use Brouwer's fixed point theorem (\comm{see e.g.,~\cite[Corollary~1.1.1]{Florenzano2003}}). That is, since $ h_{\calX} $ is a continuous map from a compact convex set $ \calX^N $ into itself, there exists a point $ x^{\sfe} \in {\calX}^N $ such that $ x^{\sfe} = h(x^\sfe) $.
  \hfill $ \Box $
\end{pf}

Next, as a tool for the convergence analysis of \eqref{eq:mpc_dynamics_conti}, we consider the entropic OT cost
\begin{equation}\label{eq:transport_cost_conti}
	\calE(x,x^\sfd) := \min_{P\in \calT(\frac{\one}{N}, \frac{\one}{N})} \ \sum_{i,j\in \bbra{N}} C_{ij} (x) P_{ij} - \varepsilon \rmH(P) \subm{,}
\end{equation}
\subm{where $ x^\sfd := [x_1^\sfd; \cdots ; x_N^\sfd] $.}
Assume that for any $ i,j\in \bbra{N} $, $ c_{T_\rmh}^i (x_i, x_j^\sfd) $ is continuously differentiable with respect to $ x_i $.
Then, thanks to the regularization, $ \calE(x,x^\sfd) $ is continuously differentiable~\cite[Eq. (9.6)]{Peyre2019} with respect to $ x $\subm{,} and
\begin{align}
	\nabla_{x_i} \calE(x, x^\sfd) = \sum_{j=1}^N P_{ij}^* (x) \nabla_{x_i} c_{T_\rmh}^i (x_i, x_j^\sfd) . \label{eq:grad_E}
\end{align}
This is in clear contrast to the case without the entropy regularization ($ \varepsilon = 0 $), in which the optimal coupling $ P^*(x) $ is not continuous similarly to the optimal permutation $ \sigma(\cdot;x) $, and thus $ \calE $ with $ \varepsilon = 0 $ is not differentiable. This difference is crucial for analyzing the \subm{global convergence property} of \eqref{eq:mpc_dynamics_conti} as shown in Theorem~\ref{thm:lasalle}.
If \comm{navigator functions} $ \{x_i^\tmp\} $ are designed appropriately, then it is expected that the state $ x(t) $ following~\eqref{eq:mpc_dynamics_conti} moves in a direction where the cost $ \calE(x(t),x^\sfd) $ decreases.
In fact, the following result shows that this is indeed the case and, as a result, ensures the convergence to the set of equilibria.
We say that $ x(t) $ converges to a set $ \calM \subset \bbR^{nN} $ as $ t\rightarrow\infty $ if for each $ \epsilon > 0 $, there exists $ \tau > 0 $ such that $\inf_{p\in \calM} \| x(t) - p \| < \epsilon $ for all $t \ge \tau $.
The proof of Theorem~\ref{thm:lasalle} is shown in Appendix~\ref{app:lasalle}.

\begin{theorem}\label{thm:lasalle}
  Suppose that Assumption~\ref{ass:positive_cost} holds.
  Assume that for any $ i\in \bbra{N} $ and $ \hat{x}\in \bbR^n $, $ c_{T_\rmh}^i (x_i, \hat{x}) $ is continuously differentiable with respect to $ x_i $ and $ T_\rmh $. Also assume that for any $ i,j\in \bbra{N} $,
  \begin{equation}\label{eq:unbounded}
    c_{T_\rmh}^i (x_i, x_j^\sfd) \rightarrow +\infty ~~ \text{as} \ \|x_i\| \rightarrow +\infty .
  \end{equation}
  Assume further that for any $ i\in \bbra{N} $, there exists a constant $ a_i > 0 $ such that for any $ x = [x_1;\cdots;x_N] \in \bbR^{nN} $,
	\begin{equation}\label{eq:tmp_condition}
		\sum_{j=1}^{N} P_{ij}^* (x) \nabla_{1} c_{T_\rmh}^i (x_i, x_j^\sfd) = a_i \nabla_{1} c_{T_\rmh}^i \left(x_i, x_i^\tmp (P^* (x))\right) ,
	\end{equation}
	where $ \nabla_1 c_{T_\rmh}^i $ denotes the gradient of $ c_{T_\rmh}^i $ with respect to the first variable.
	Then, for any initial state $ x(0) = x^0 \in \bbR^{nN} $, \blu{the solution} $ x(t) $ \blu{of \eqref{eq:mpc_dynamics_conti}} converges to the set $ \calR = \left\{x^\sfe \in \bbR^{nN} : x_i^\sfe  = x_i^\tmp (P^* (x^\sfe)), \forall i \in \bbra{N} \right\} $ as $ t\rightarrow \infty $.
	\hfill $ \diamondsuit $
\end{theorem}

\begin{rmk}
  In the proof of Theorem~\ref{thm:lasalle}, the linearity of the system~\eqref{eq:linear_dynamics_conti} is not used at all.
  Hence, the same proof works for nonlinear systems of the form:
  \[
    \dot{x}_i = f_i (x_i) + g_i(x_i) u_i ,
  \]
  where $ g_i(x_i) $ is square and invertible for any $ x_i\in \bbR^n $.
  \hfill $ \diamondsuit $
\end{rmk}
\begin{rmk}
  Assume that $\calR $ consists of only isolated points. Then by Theorem~\ref{thm:lasalle}, for any initial state, $ x(t) $ converges to one of the equilibrium points in $ \calR $ as $ t\rightarrow \infty $.
  \hfill $ \diamondsuit $
\end{rmk}
\begin{rmk}\label{rmk:nonuniform}
	Let us consider the general case where the number of agents $ N $ is not equal to the number of targets $ M $, and the agents and targets have mass distributions that are not necessarily uniform, i.e., $ \bfa \neq \bfone_N/N, \ \bfb \neq \bfone_M/M $. Then, $ P^* $ in the dynamics~\eqref{eq:mpc_dynamics_conti} is replaced by
	\begin{equation}
	\bar{P}^* (x) := \argmin_{P\in \calT(\bfa, \bfb)} \ \sum_{i\in \bbra{N}, j\in \bbra{M}} C_{ij} (x) P_{ij} - \varepsilon \rmH(P) .
	\end{equation}
	Even in this case, the same proof as in Theorem~\ref{thm:lasalle} works. That is, under the assumptions in Theorem~\ref{thm:lasalle}, for any initial state $ x(0) = x^0 \in \bbR^{nN} $, $ x(t) $ converges to $ \calR $ as $ t\rightarrow \infty $.
	In other words, Sinkhorn MPC can be applied to general OT problems whereas MPC with e.g., the Hungarian algorithm works only for OT problems that are equivalent to assignment problems. This is one of the advantages of the proposed method.
	\hfill $ \diamondsuit $
\end{rmk}

The condition \eqref{eq:tmp_condition} gives a guideline for the design of $ x_i^\tmp $.
However, it is not necessary for $ x_i^\tmp $ to satisfy \eqref{eq:tmp_condition} to ensure the convergence. In fact, in Section~\ref{sec:example}, we observed that for the non-quadratic cost~\eqref{eq:l1_cost}, the barycentric projection~\eqref{eq:barycentric_target}, which is not desinged based on \eqref{eq:tmp_condition} achieves the transport to the target distribution.

Next, as a specific example, we consider a quadratic cost
\begin{equation}\label{eq:quadratic}
	\ell_i (x_i, u_i; x_j^\sfd) = \|u_i + B_i^{-1} A_i x_j^\sfd\|^2 .
\end{equation}
Then the corresponding control law $ u_i^\mpc $ and the cost function $ c_{T_\rmh}^i $ can be written as follows~\cite[Section~3.3,~pp.~138-140]{Lewis2012}:
\begin{align}
  &u_i^{\rm MPC} (x_i, \red{\hat{x}}) = - B_i^\top \scrG_i (x_i - \red{\hat{x}} ) - B_i^{-1} A_i \hat{x}, \label{eq:mpc_quad_conti}\\
    &c_{T_\rmh}^i (x_i, \hat{x}) = \| x_i - \hat{x} \|_{\scrG_i}^2, ~~ \forall x_i, \hat{x} \in \bbR^n,\ \forall i \in \bbra{N}, \label{eq:quad_value_func_conti}
  \end{align}
  where
  \begin{align*}
  & \scrG_i := \Ee^{A_i^\top T_\rmh} G_{i,T_\rmh}^{-1} \Ee^{A_i T_\rmh}, \ G_{i,T_\rmh} := \int_0^{T_\rmh} \Ee^{A_i t} B_i B_i^\top  \Ee^{A_i^\top t} \rmd t .
  \end{align*}
  Thus, the condition \eqref{eq:tmp_condition} is rewritten as
  \begin{align*}
    \frac{2}{N} \scrG_i \biggl(x_i - N\sum_j P_{ij}^* (x) x_j^\sfd \biggr) = 2a_i \scrG_i \left(x_i - x_i^\tmp (P^*(x))\right) .
  \end{align*}
  This is fulfilled by $ a_i = 1/N $ and the barycentric projection~\eqref{eq:barycentric_target}.

The cost function \eqref{eq:quadratic} does not satisfy Assumption~\ref{ass:positive_cost} because it does not depend on the variable $ x_i $. Nevertheless, similarly to Theorem~\ref{thm:lasalle}, the following convergence result holds. The proof is given in Appendix~\ref{app:lasalle}.
\begin{corollary}\label{cor:convergence}
  Consider the quadratic cost \eqref{eq:quadratic} with $ {\bbX_i = \bbU_i = \bbR^n} $ for all $ i\in \bbra{N} $ and the barycentric target~\eqref{eq:barycentric_target}.
  Then, for any initial state $ x(0) = x^0 \in \bbR^{nN} $, \blu{the solution} $ x(t) $ \blu{of \eqref{eq:mpc_dynamics_conti}} converges to the set \blu{$ \calR = \{x^\sfe \in \bbR^{nN} : x_i^\sfe  = N \sum_{j=1}^N P_{ij}^*(x^\sfe) x_j^\sfd  , \forall i \in \bbra{N}\} $} as $ t\rightarrow \infty $.
	\hfill $ \diamondsuit $
\end{corollary}

The above result justifies that the barycentric projection~\eqref{eq:barycentric_target}, which \subm{is} typically used to approximate Monge's OT map from a coupling matrix, gives an appropriate direction where the cost $ \calE(x(t), x^\sfd) $ becomes smaller under Sinkhorn MPC for the quadratic cost~\eqref{eq:quadratic}.

\res{Theorem~\ref{thm:lasalle} and Corollary~\ref{cor:convergence} ensure the global convergence of the proposed method that uses the optimal coupling $ P^*(x(t)) $ at each time.
Hence, in terms of the convergence, it is desirable to perform a sufficiently large number of Sinkhorn iterations at each time to obtain a coupling close enough to $ P^*(x(t)) $.
On the other hand, as observed in Fig.~\ref{fig:iteration}, once we perform a sufficiently large number of iterations at some point, we can obtain a coupling close enough to $ P^*(x(t)) $ at later times by a smaller number of iterations. 
This implies that we can reduce the computational burden at later times while maintaining the convergence property.}

\section{Boundedness and local asymptotic stability for Sinkhorn MPC with a quadratic cost}\label{sec:linear_system}
In this section, we consider the general case where the number of the Sinkhorn iterations $ S $ is not necessarily large, and thus $ P[k] $ may not be close enough to the optimal coupling $ P^* (\bmx[k]) $.
Then, we elucidate the fundamental properties of Sinkhorn MPC on discrete-time systems~\subm{\eqref{eq:linear_dynamics_disc}}.
Specifically, we reveal that even when $ S $ is small, the ultimate boundedness and the local asymptotic stability for Sinkhorn MPC hold for \modi{the quadratic cost \eqref{eq:energy}} and $ \bbX_i = \bbU_i = \bbR^n $. \rev{Hereafter, we assume the invertibility of $ \bmB_i $.}
In addition, for notational simplicity, we deal only with the case where just one Sinkhorn iteration is performed at each time, i.e., $ S[k] = 1, \forall k $. Nevertheless, by similar argument, all \rev{of} the results in this section are still valid when more iterations are performed.
For convenience, we recall the dynamics under Sinkhorn MPC with $ S \subm{[k]}  \equiv 1 $:
\begin{align}
	&\bmx_i[k+1] = \bar{\bmA}_i \bmx_i [k] +  (I - \bar{\bmA}_i) x_i^\tmp (P[k]),\ \forall i \in \bbra{N}, \label{eq:SMPC_linear_x_1} \\
	&P[k] = \alpha[k+1]^\diagbox \bmK(\bmx[k]) \beta[k]^\diagbox, \label{eq:SMPC_P_1} \\
	&\alpha[k+1] = \bfone_N/N \oslash \left[\bmK(\bmx[k])\beta[k] \right], \label{eq:SMPC_alpha_1} \\
	&\beta[k] = \bfone_N/N \oslash \left[\bmK(\bmx[k])^\top \alpha[k] \right], \label{eq:SMPC_beta_1}   \\
	&\bmx_i [0] = x_i^0, \ \alpha[0] = \alpha_0 . \nonumber 
\end{align}

\subsection{Ultimate boundedness for Sinkhorn MPC}\label{subsec:bounded}
Here, we assume that for the codomain $ \calX $ of $ x_i^\tmp $, there exists a constant $ \bar{r} > 0 $ such that
\begin{equation}\label{ineq:x_upp}
	\|x \| \le  \bar{r} , \ \forall x \in \calX .
\end{equation}
For example, if $ \calX $ is the convex hull of $ \{x_j^\sfd\}_j $, we can take $ \bar{r} = \max_j \|x_j^\sfd\| $.
It is known that, under the assumption that $ \bmB_i $ is invertible, $ \bar{\bmA}_i $ is stable, i.e., \blue{the spectral radius $ \rho_i $ of $ \bar{\bmA}_i $ satisfies $ \rho_i < 1 $} \cite[Corollary~1]{Kwon1975}.
Using this fact, we derive the ultimate boundedness of \eqref{eq:SMPC_linear_x_1} with \eqref{eq:SMPC_P_1}--\eqref{eq:SMPC_beta_1}.

\subm{
\begin{proposition}\label{prop:bounded}
	Assume that there exists a constant $ \bar{r} > 0 $ satisfying \eqref{ineq:x_upp}.
	Then, for any $ \delta > 0, \{x_i^0\}_i $, and $ \{\nu_i\}_i $ satisfying $ \nu_i > 0, \rho_i + \nu_i < 1, \forall i\in \bbra{N} $, there exist $ \kappa_i(\nu_i) > 0, i\in \bbra{N} $ and $ \tau(\delta, \{x_i^0\},\{\nu_i \} ) \in \bbZ_{>0} $ such that the solution $ \{\bmx_i\}_i $ of \eqref{eq:SMPC_linear_x_1} with \eqref{eq:SMPC_P_1}--\eqref{eq:SMPC_beta_1} satisfies
	\begin{equation}\label{ineq:ultimate_bound}
		\|\bmx_i [k] \| < \delta  +  \frac{ \kappa_i \bar{r} \|I- \bar{\bmA}_i\|_2}{1 - (\rho_i + \nu_i ) }  , \ \forall k \ge \blue{\tau}, \ \forall i \in \bbra{N}.
	\end{equation}
\end{proposition}
\begin{pf}
	Let $ \tilde{\bmu}_i [k] := ( I - \bar{\bmA}_i) x_i^\tmp (P[k]) $. 
	Then, it follows from \eqref{ineq:x_upp} that
	\begin{align*}
	\| \tilde{\bmu}_i [k] \| \le \bar{r} \| I -  \bar{\bmA}_i \|_2, \ \forall k\in \bbZ_{\ge 0} .
	\end{align*}
	By \cite[Lemma~7.3.2]{Golub2013}, for any $ \nu_i > 0 $, there exists $ \kappa_i (\nu_i) > 0 $ such that
	\[
		\| \bar{\bmA}_i^k \|_{\subm{2}} \le \kappa_i (\rho_i + \nu_i )^k, \ \forall k \in \bbZ_{\ge 0} . 
	\]
	Hence, the desired result is straightforward from
	\begin{align*}
	\| \bmx_i [k] \| \le \| \bar{\bmA}_i^k \|_2 \| x_i^0 \| + \sum_{s = 1}^{k} \| \bar{\bmA}_i^{s-1} \|_2 \| \tilde{\bmu}_i [k-s]\| .
	\end{align*}
	\hfill $ \Box $
\end{pf}
}

\rev{We emphasize that Proposition~\ref{prop:bounded} holds for any \comm{navigator function} $ x_i^\tmp $ \fin{whose codomain $ \calX $ satisfies \eqref{ineq:x_upp}.}}

\subsection{Existence of the equilibrium points}\label{subsec:equilibrium}
In the \rev{remainder} of this section, we focus on the barycentric target~\eqref{eq:barycentric_target}.
For $ (x,\beta) \in \bbR^{nN} \times \bbR_{>0}^N $ and $ X^\sfd := [x_1^\sfd \ \cdots \ x_N^\sfd] \in \bbR^{n\times N} $, define
\begin{align*}
&f_1 (x,\beta) := \{\bar{\bmA}_i\}_i^\diagbox x \nonumber\\
&\hspace{1.6cm}+ N  \{I_n - \bar{\bmA}_i\}_i^\diagbox  (X^\sfd)^{\diagbox,N} {\rm vec}(\tilde{P}(x,\beta)), \\
&f_2 (x,\beta) := \bfone_N/N \oslash \left[ \bmK(f_1 (x,\beta))^\top (\bfone_N/N \oslash [\bmK(x) \beta]) \right], \\
&\tilde{P}(x,\beta) :=  \left( \bfone_N/N \oslash [\bmK(x)\beta]  \right)^\diagbox \bmK(x)  \beta^\diagbox. 
\end{align*}
Then, the collective dynamics \eqref{eq:SMPC_linear_x_1}--\eqref{eq:SMPC_beta_1} is \blu{rewritten as}
\begin{align}
&\bmx[k+1] = f_1 (\bmx[k], \beta[k])  , \label{eq:x_global}\\
&\beta[k+1] = f_2(\bmx[k], \beta[k]) . \label{eq:beta_global}
\end{align}

A point $ x^{\sfe} = \rev{[x_1^\sfe; \ \cdots \ ;x_N^\sfe ]} \in \bbR^{nN} $ is an equilibrium of \eqref{eq:x_global},~\eqref{eq:beta_global} if and only if
\begin{align}
&\hspace{-0.2cm}(I_n - \bar{\bmA}_i) \biggl( x_{i}^{\sfe} - N \sum_{j=1}^{N} \bmP_{ij}^* (x^\sfe) x_j^\sfd \biggr) = 0 , \ \forall i \in \bbra{N}. \nonumber
\end{align}
Here note that $ \bmP^* $ satisfies
\begin{align}
&\hspace{-0.2cm} \bmP_{ij}^* (x) =  \alpha_{i}^* \bmK_{ij}(x) \beta_{j}^*,\ \alpha^*, \beta^* \in \bbR_{>0}^N,\\
&\hspace{-0.2cm}\alpha^* = \bfone_N/N \oslash \left[ \bmK(x) \beta^*	\right],  \beta^* = \bfone_N/N \oslash \left[ \bmK(x)^\top \alpha^*	\right]. \label{eq:optimal_beta}
\end{align}
The stability of $ \bar{\bmA}_i $ implies that it has no eigenvalue \rev{equal to} $ 1 $, and therefore $ I_n - \bar{\bmA}_i $ is invertible. Thus, the necessary and sufficient condition for the equilibria is given by
\begin{equation}\label{eq:equilibrium_iff}
x_{i}^{\sfe} - N\sum_{j=1}^{N} \bmP_{ij}^* (x^\sfe) x_j^{\sfd} = 0, \ \forall i\in \bbra{N} ,
\end{equation}
which coincides with \eqref{eq:equi_general}.
Similarly to Proposition~\ref{prop:fixed_point_general}, we show the existence of an equilibrium.
\begin{corollary}\label{prop:fixed_point}
	The dynamics \eqref{eq:x_global}, \eqref{eq:beta_global} has at least one equilibrium point $ (x^\sfe, \beta^\sfe) \in \bbR^{nN} \times (\bbR_{> 0}^N / {\sim}) $.
\end{corollary}
\begin{pf}
	Note that if a point $ x^{\sfe}\in \bbR^{nN} $ satisfies \eqref{eq:equilibrium_iff}, the corresponding $ \beta^\sfe \in \bbR_{> 0}^N / {\sim} $ is uniquely determined by $ {\beta^\sfe = \beta^* }$ in \eqref{eq:optimal_beta} with $ x= x^\sfe $~\cite[Theorem~4.2]{Peyre2019}.
	Note also that for any $ i \in \bbra{N}$ and any $ x \in \bbR^{nN} $, $ N\sum_{j=1}^N \bmP_{ij}^* (x) x_j^\sfd $ belongs to the convex hull $ \calX $ of $ \{x_j^\sfd\}_j $.
	Then, by the same argument as in the proof of Proposition~\ref{prop:fixed_point_general}, we obtain the desired result.
	\hfill $ \Box $
\end{pf}

Sometimes, in order to emphasize the dependence of $ (x^\sfe, \beta^\sfe) $ on $ \varepsilon $, we write $ (x^\sfe (\varepsilon), \beta^\sfe (\varepsilon)) $.

\subsection{Local asymptotic stability for Sinkhorn MPC}
Next, we analyze the stability of the equilibrium points.
For this purpose, the following lemma is crucial when $ \varepsilon $ is small.
The proof is shown in Appendix~\ref{app:convergence}.
\begin{lemma}\label{lemma:exponential_convergence}
	Assume that $ x_i^\sfd \neq x_j^\sfd $ for all $(i,j), \ i\neq j $, \blu{and $ \bmA_i $ is invertible for all $ i\in \bbra{N} $.}
	For a permutation $ \sigma : \bbra{N} \rightarrow \bbra{N} $, define $ x^{\sfd} (\sigma) := \rev{[x_{\sigma(1)}^{\sfd}; \ \cdots \ ;x_{\sigma(N)}^{\sfd} ]} $ and a permutation matrix $ P^\sigma = (P_{ij}^\sigma)$ as $ P_{ij}^\sigma := 1/N $ if $ j= \sigma(i) $, and $ 0 $, otherwise.
	Then \subm{for any permutation $ \sigma $,} there exists an equilibrium $ (x^\sfe(\varepsilon), \beta^\sfe(\varepsilon)) $ of \eqref{eq:x_global}, \eqref{eq:beta_global} such that $ x^\sfe (\varepsilon) $ and $ P^* (x^\sfe (\varepsilon)) $ converge exponentially to $ x^{\sfd} (\sigma) $ and $ P^\sigma $, respectively, as $ \varepsilon \rightarrow +0 $, i.e., there exists $ \zeta > 0 $ such that
	\[
		\lim_{\varepsilon \rightarrow +0} \frac{\|\eta(\varepsilon)  \|_2}{\exp (-\zeta/\varepsilon)} = 0
	\]
	for $ \eta (\varepsilon) = x^\sfe (\varepsilon) - x^\sfd (\sigma) $ and $ \eta (\varepsilon) = P^*(x^\sfe (\varepsilon)) - P^\sigma$.
	\hfill $ \diamondsuit $
\end{lemma}

Denote by $ {\rm Exp}(\sigma) $ the set of all equilibria \fin{$ (x^\sfe (\cdot), \beta^\sfe (\cdot)) $} of \eqref{eq:x_global}, \eqref{eq:beta_global} having the property in Lemma~\ref{lemma:exponential_convergence} for a permutation $ \sigma $.

For $ \bar{P} \in \bbR^{N\times N} $ and $ x  = \rev{[x_1; \ \cdots \ ;x_N]} \in \bbR^{nN} $, define
\begin{align*}
	V_{\subm{\bar{P}}}(x) := \sum_{i=1}^{N} \Bigl\| x_i - N \sum_{j=1}^{N} \bar{P}_{ ij} x_j^\sfd \Bigr\|_{\calG_i}^2 .  
\end{align*}
Then, $ V_{\subm{\bar{P}}} $ is a Lyapunov function of \eqref{eq:x_global} where $ \tilde{P}(x,\beta) $ is fixed by $ \bar{P} $~\cite{Mayne2000}.
Indeed, we have
\begin{align*}
	&V_{\subm{\bar{P}}} (\bmx[k+1]) - V_{\subm{\bar{P}}} (\bmx[k]) \le - \sum_{i=1}^N W_{1,i} (\bmx_i[k], \bar{P}) \\
	&W_{1,i} (x_i, \bar{P}) := \Bigl\| \bmB_i^\top (\bmA_i^\top)^{\tau_\rmh -1} \bmG_{i,\tau_\rmh}^{-1} \bmA_i^{\tau_\rmh} \\
	&\hspace{4cm} \times \bigl(x_i - N\sum_{j=1}^N \bar{P}_{ij} x_j^\sfd \bigr) \Bigr\|^2 . 
\end{align*}
Given an equilibrium $ (x^\sfe, \beta^\sfe) $, let us take the optimal coupling $ P^\sfe := \bmP^* (x^{\sfe}) $ as $ \bar{P} $, and for $ \gamma > 0 $, define
\begin{equation}\label{eq:V}
V(x, \beta) :=  V_{\subm{P^\sfe}} (x) +  \gamma \dhil ( \beta, \beta^\sfe ), \ (x,\beta)\in \bbR^{nN} \times  (\bbR_{>0}^N / {\sim} ) . 
\end{equation}
The following theorem follows from the fact that, for sufficiently small or large $ \varepsilon > 0$ and large $ \gamma > 0 $, $ V $ behaves as a Lyapunov function of \eqref{eq:x_global},~\eqref{eq:beta_global} with respect to $ (x^\sfe, \beta^\sfe) $.
We give the proof in Appendix~\ref{app:convergence}.

\begin{theorem}\label{thm:stability}
	Assume that for all $ i\in \bbra{N} $, $ \bmA_i $ is invertible.
	Then the following hold:
	\begin{itemize}
		\item[(i)] Assume that $ (x^\sfe, \beta^\sfe) $ is an isolated equilibrium of \eqref{eq:x_global}, \eqref{eq:beta_global}. Then, for a sufficiently large $ \varepsilon > 0 $, $ (x^\sfe, \beta^\sfe) $ is locally asymptotically stable.
		\item[(ii)]  
		Assume that $ x_i^\sfd \neq x_j^\sfd $ for all $(i,j), \ i\neq j $. \fin{Assume further that for some $ \varepsilon' > 0 $, $ (x^\sfe (\varepsilon'), \beta^\sfe(\varepsilon')) $ is an isolated equilibrium of \eqref{eq:x_global}, \eqref{eq:beta_global} and for some permutation $ \sigma $, $ (x^\sfe (\cdot), \beta^\sfe(\cdot)) \in {\rm Exp}(\sigma) $.}
		Then, for sufficiently small $ \varepsilon >0$, $ (x^\sfe (\varepsilon), \beta^\sfe(\varepsilon)) $ is locally asymptotically stable. \hfill$ \diamondsuit $
	\end{itemize}
\end{theorem}

\subsection{Interpretation of Sinkhorn MPC as an alternating descent and ascent method}
Lastly, we give an interpretation of Sinkhorn MPC.
We continue to work on the quadratic cost~\eqref{eq:energy} with the barycentric target~\eqref{eq:barycentric_target}.
	First, similarly to \eqref{eq:transport_cost_conti} for the continuous-time systems, consider the entropic OT cost
	\begin{equation}\label{eq:OT_cost_disc}
		E(x,x^\sfd) := \min_{P\in \calT(\frac{\one}{N}, \frac{\one}{N})} \ \sum_{i,j\in\bbra{N}} \bmC_{ij} (x) P_{ij} - \varepsilon \rmH(P),
	\end{equation}
	where $ \bmC_{ij}(x) = \|x_i - x_j^\sfd \|_{\calG_i}^2 $.
	Since it holds
	\begin{align*}
		\nabla_{x_i} E(x, x^\sfd) = 2\calG_i \sum_{j=1}^{N} \bmP_{ij}^* (x) (x_i - x_j^\sfd) ,
	\end{align*}
	\blu{a point $ x $ satisfying the following is a stationary point of $ E $.}
	\begin{equation}\label{eq:stationary_point}
		x_i - N \sum_{j=1}^{N} \bmP_{ij}^* (x) x_j^\sfd = 0, \ \forall i\in \bbra{N} .
	\end{equation}
	This is exactly the condition \eqref{eq:equilibrium_iff} for the equilibrium points $ x^{\sfe} $ of \eqref{eq:x_global}, \eqref{eq:beta_global}.
	\subm{Hence, Sinkhorn MPC can be viewed as a cost-effective search method to find the stationary points of the associated entropic OT cost $ E $.}

	Next, we introduce the dual problem associated \subm{with} \eqref{eq:OT_cost_disc}:
	\begin{align*}
		& E(x,x^\sfd) = \max_{{\bf f},{\bf g}\in \bbR^N} \ Q({\bf f},{\bf g}; x), \\
		&Q({\bf f},{\bf g}; x) := {\bf f}^\top (\bfone_N / N) + {\bf g}^\top (\bfone_N / N) - \varepsilon (\Ee^{{\bf f}/\varepsilon})^\top \bmK(x) \Ee^{{\bf g}/\varepsilon} ,
	\end{align*}
	where $ \Ee^{{\bf f}/\varepsilon} \in \bbR^{N} $ denotes the element-wise exponential of $ {\bf f}/\varepsilon $.
	Let $ {\bf f}[k] = \varepsilon \log (\alpha[k]), {\bf g}[k] = \varepsilon \log (\beta[k]) $ for scaling variables of the Sinkhorn algorithm~\eqref{eq:Sinkhorn_algorithm} where $ (\log(\alpha))_i = \log (\alpha_i) $. Then, the Sinkhorn iterations \eqref{eq:Sinkhorn_algorithm} are equivalent to a block coordinate ascent~\cite[Remark~4.21]{Peyre2019}, which updates alternatively $ {\bf f} $ and $ {\bf g} $ to cancel the respective gradients
	\begin{align*}
		&\nabla_{\bf f} Q({\bf f},{\bf g}; x) = \bfone_N/N - \Ee^{{\bf f}/\varepsilon} \odot (\bmK(x)\Ee^{{\bf g}/\varepsilon}), \\
		&\nabla_{\bf g} Q({\bf f},{\bf g}; x) = \bfone_N/N - \Ee^{{\bf g}/\varepsilon} \odot (\bmK(x)^\top \Ee^{{\bf f}/\varepsilon}) ,
	\end{align*}
	where $ \odot $ denotes element-wise multiplication.
	On the other hand, the gradient with respect to $ x_i $ is
	\begin{align}
		&\nabla_{x_i} Q({\bf f}[k+1], {\bf g}[k] ; x) \nonumber\\
		&= \sum_{j\in \bbra{N}} \alpha_i[k+1] \beta_j[k] \exp\left( - \frac{\bmc_{\tau_\rmh}^i (x_i,x_j^\sfd)}{\varepsilon}  \right) \nabla_{x_i} \bmc_{\tau_\rmh}^i (x_i,x_j^\sfd) \nonumber\\
		&= \sum_{j\in \bbra{N}} P_{ij} [k] \nabla_{x_i} \bmc_{\tau_\rmh}^i (x_i, x_j^\sfd),
	\end{align}
	which has the same form as \eqref{eq:grad_E}. Now, let us consider the case when the dynamics~\eqref{eq:linear_dynamics_conti} is well approximated by the discretized system~\eqref{eq:linear_dynamics_disc}. 
	Then using the same derivation as for \eqref{eq:lyap_decrease} in the proof of Theorem~\ref{thm:lasalle}, along the trajectory $ x(t) $ following
	\begin{equation}
		\dot{x}_i (t) = A_i x_i (t) + B_i u_i^{\rm MPC} \left(x_i(t), x^{\rm tmp}_{i} (P[k]) \right) ,
	\end{equation}
	we have
	\begin{align}
		&\frac{\rmd}{\rmd t} \calQ\left(\bff[k+1], \bfg[k] ; x(t)\right) \nonumber\\
		&\qquad\begin{cases}
			< 0, & x_i(t)  \neq x_i^\tmp (\subm{P[k]}), \ \exists i\in \bbra{N}, \\
			= 0, & x_i(t)  = x_i^\tmp (\subm{P[k]}), \ \forall i\in \bbra{N} ,
		\end{cases}
	\end{align}
	where $ \calQ $ is the continuous-time version of $ Q $:
	\begin{align*}
		\calQ({\bf f},{\bf g}; x) &:= {\bf f}^\top (\bfone_N / N) + {\bf g}^\top (\bfone_N / N) - \varepsilon (\Ee^{{\bf f}/\varepsilon})^\top K(x) \Ee^{{\bf g}/\varepsilon} \\
		&\simeq Q({\bf f},{\bf g}; x).
	\end{align*}
	Therefore, if the sampling time is small, the state trajectory $ \bmx[k] $ moves in a direction where $ Q(\bff[k+1], \bfg[k]; \bmx[k]) $ decreases.
	Note that the above argument applies to the general cost under the assumptions in Theorem~\ref{thm:lasalle}.

	In summary, Sinkhorn MPC can be interpreted as an alternating descent and ascent method to seek a solution of the minimax problem
	\[
		\min_{x\in \bbR^{nN}} \max_{{\bf f},{\bf g}\in \bbR^N} \ Q({\bf f},{\bf g};x) ,
	\]
	where the minimizers satisfy \eqref{eq:stationary_point}.

\section{Conclusion}\label{sec:conclusion}
\magenta{
In this paper, we presented the concept of Sinkhorn MPC, which integrates MPC and the Sinkhorn algorithm to achieve scalable, cost-effective transport over dynamical systems.
The numerical examples described the usefulness of the proposed method.
Moreover, thanks to the entropy regularization, under some assumptions, we ensured the global convergence for Sinkhorn MPC, which is one of the most important properties of transport algorithms.
Furthermore, for linear systems with a quadratic cost, we analyzed the ultimate boundedness and the local asymptotic stability for Sinkhorn MPC based on the stability of the constrained MPC and the conventional Sinkhorn algorithm.

On the other hand, in the numerical example, we observed that the regularization parameter plays a key role in the trade-off between the stationary and transient behaviors for Sinkhorn MPC.
Hence, an important direction for future work is to investigate the design of a time-varying regularization parameter to balance the trade-off.
In addition, although we focused on the case where an OT problem is equivalent to an assignment problem, in Remark~\ref{rmk:nonuniform}, we mentioned that the convergence result for Sinkhorn MPC still holds in more general settings. Hence, it is also interesting to explore applications of Sinkhorn MPC for general OT problems.
\res{In this paper, for simplicity, we assumed the invertibility of $ B_i $ for the convergence analysis. Possible relaxation of this assumption will be reported in a future publication.}

}

\begin{ack}
	This work was supported in part by JSPS KAKENHI Grant Numbers JP21J14577, JP21H04875, JST, ACT-X Grant Number JPMJAX2102, and the joint project of Kyoto University and Toyota Motor Corporation, titled ``Advanced Mathematical Science for Mobility Society.''
\end{ack}

\appendix

\section{Proofs of Theorem~\ref{thm:lasalle} and Corollary~\ref{cor:convergence}}\label{app:lasalle}
\subsection{Proof of Theorem~\ref{thm:lasalle}}
We prove Theorem~\ref{thm:lasalle} by using LaSalle's invariance principle~\cite[Theorem~4.4]{Khalil2002}.
	\begin{proposition}\label{prop:invariance}
		Suppose that there exist a constant $ d \in \bbR $ and a continuously differentiable function $ V : \bbR^{nN} \rightarrow \bbR $ such that the sublevel set $ \Omega_V (d) := \{ x \in \bbR^{nN} : V(x) \le d\} $ is bounded, and $ \frac{\rmd }{\rmd t} V(x(t)) |_{x(t) = x'} \le 0 $ for all $ x' \in \Omega_V (d) $. Let 
		\[
			\calR_V (d) := \left\{x' \in \Omega_V(d) : \frac{\rmd V(x(t))}{\rmd t} \biggl|_{x(t) = x'}  = 0 \right\},
		\]
		and let $ \calM $ be the largest invariant set in $ \calR_V (d) $. Then every solution of \eqref{eq:mpc_dynamics_conti} starting in $ \Omega_V (d) $ converges to $ \calM $ as $ t\rightarrow \infty $.
		\hfill $ \diamondsuit $ 
	\end{proposition}

	As a candidate for the above function $ V $, we choose $ \calE(x,x^\sfd) $.
  The time derivative of $ \calE (x(t),x^\sfd) $ along the trajectory of \eqref{eq:mpc_dynamics_conti} is given by
  \begin{align}
    &\frac{\rmd}{\rmd t} \calE (x(t), x^\sfd) = \nabla_x \calE (x(t), x^\sfd)^\top \dot{x} (t) \nonumber\\
    &= [\nabla_{x_1} \calE(x(t),x^\sfd);\cdots ; \nabla_{x_N} \calE(x(t), x^\sfd) ]^\top \dot{x} (t) \nonumber\\
    &= \sum_{i=1}^N \sum_{j=1}^N P_{ij}^* (x(t)) \nabla_{x_i} c_{T_\rmh}^i (x_i (t), x_j^\sfd)^\top \nonumber\\
    &\quad \times \left( A_i x_i(t) + B_i u_i^{\rm MPC} \Bigl(x_i(t), x_i^\tmp \bigl(P^* (x(t))\bigr)\Bigr) \right) ,
  \end{align}
  where we used \eqref{eq:grad_E}.
	By the same argument as in the proof of \cite[Theorem~1]{Chen1982}, which derives the stability for MPC with a terminal equality constraint, under the differentiability of $ c_{T_\rmh}^i(x_i,\hat{x}) $ with respect to $ x_i $ and $ T_\rmh $, it can be shown that
\begin{align}
	&\nabla_{x_i} c_{T_\rmh}^i (x_i,\hat{x})^\top \left(A_i x_i  + B_i u_i^{\rm MPC} (x_i,\hat{x}) \right) \nonumber\\
	&\le  - \ell_i \left(x_i, u_i^{\rm MPC} (x_i,\hat{x}); \hat{x} \right) , \ \forall i \in \bbra{N},\  \forall x_i, \hat{x} \in \bbR^n . \label{eq:cost_bound}
\end{align}
By \eqref{eq:tmp_condition} and \eqref{eq:cost_bound}, it holds
\begin{align*}
	&\frac{\rmd }{\rmd t} \calE (x(t), x^\sfd) = \sum_{\blu{i}} a_i \nabla_{1} c_{T_\rmh}^i \left(x_i(t), x_i^\tmp \subm{\bigl(}P^* (x(t))\bigr)\right)^\top \\
  &\times  \left( A_i x_i(t) + B_i u_i^{\rm MPC} \Bigl(x_i(t), x_i^\tmp \subm{\bigl(}P^* (x(t))\bigr)\Bigr) \right) \\
	&\le -  \sum_{\blu{i}} a_i \ell_i \Bigl(x_i(t), u_i^{\rm MPC} \bigl(x_i(t), x_i^\tmp \subm{\bigl(}P^* (x(t))\bigr) \bigr) \\
  &\qquad\qquad\qquad ; x_i^\tmp \subm{\bigl(}P^* (x(t))\bigr) \Bigr) .
\end{align*}
Therefore, by \eqref{eq:positive_cost} in Assumption~\ref{ass:positive_cost},
\begin{align}\label{eq:lyap_decrease}
	\frac{\rmd }{\rmd t} \calE (x(t), x^\sfd)
	\begin{cases}
		< 0, & x_i(t)  \neq x_i^\tmp \subm{\bigl(}P^* (x(t))\bigr), \exists i\in \bbra{N}, \\
		= 0, & x_i(t)  = x_i^\tmp \subm{\bigl(}P^* (x(t))\bigr), \forall i\in \bbra{N} .
	\end{cases}
\end{align}

Next, we show that for any $ d\in \bbR $ such that the \blu{sublevel} set $ \Omega_\calE (d) := \{x\in \bbR^{nN} : \calE(x,x^\sfd) \le d\} $ is non-empty, $ \Omega_\calE (d) $ is bounded.
\subm{Since for all \las{$ P \in \calT(\bfone_N/N, \bfone_N/N) $},
\[
	\rmH(P) \le -\biggl(\sum_{i,j} \frac{1}{N^2} \log \biggl(\frac{1}{N^2} \biggr) \biggr) + 1 = 2\log N + 1 ,
\]
}it holds for any $ x\in \bbR^{nN} $,
\begin{align*}
	&\calE (x,x^\sfd) \ge \sum_{i,j\in \bbra{N}} P_{ij}^* (x)  c_{T_\rmh}^i (x_i,x_j^\sfd) - \varepsilon (2\log N + 1) \nonumber\\
	&\ge \sum_{i\in \bbra{N}} \frac{1}{N} \min_{j} c_{T_\rmh}^i (x_i, x_j^\sfd) - \varepsilon (2\log N + 1) =: \underline{\calE} (x,x^\sfd) .
\end{align*}
Hence, for any $ d \in \bbR $,
\begin{align*}
	\Omega_\calE (d) \subseteq \{ x\in \bbR^{nN} : \underline{\calE} (x,x^\sfd) \le d \} =: \Omega_{\underline{\calE}} (d) .
\end{align*}
In addition, by \eqref{eq:unbounded}, $ \Omega_{\underline{\calE}} (d) $ is bounded, and therefore $ \Omega_\calE (d) $ is also bounded.

For any $ d\in \bbR $, let
\begin{align*}
  \calR_\calE (d) &:= \left\{x' \in \Omega_\calE (d) : \frac{\rmd \calE(x(t),x^\sfd)}{\rmd t} \biggl|_{x(t) = x'}  = 0 \right\} \\
  &= \left\{x'\in \Omega_\calE (d) : x'_i = x_i^\tmp (P^* (x')), \ \forall i \in \bbra{N} \right\} . 
\end{align*}
Since any point in $ \calR_\calE (d) $ is an equilibrium of \eqref{eq:mpc_dynamics_conti} by \eqref{eq:equi_general}, the largest invariant set in $ \calR_\calE (d) $ is $ \calR_\calE (d) $ itself.
Therefore, by \eqref{eq:lyap_decrease} and Proposition~\ref{prop:invariance}, for any $ x(0) = x^0 \in \Omega_\calE (d) $, $ x(t) $ converges to the largest invariant set $ \calR_\calE (d) $. By the arbitrariness of $ d $, we obtain the desired result.

\subsection{Proof of Corollary~\ref{cor:convergence}}
Note that $ c_{T_\rmh}^i (x_i,\hat{x}) $ given by \eqref{eq:quad_value_func_conti} is continuously differentiable with respect to $ x_i $ and $ T_\rmh $, and satisfies \eqref{eq:unbounded}. Then by the same argument as in the proof of Theorem~\ref{thm:lasalle}, we obtain
\begin{align*}
	\frac{\rmd }{\rmd t} \calE (x(t), x^\sfd) &\le - \frac{1}{N} \sum_{\blu{i}} \Bigl\| u_i^\mpc \Bigl(x_i(t),x_i^\tmp \bigl(P^*(x(t))\bigr) \Bigr) \\
	&\qquad + B_i^{-1} A_i x_i^\tmp \bigl(P^*(x(t))\bigr) \Bigr\|^2 .
\end{align*}
Let $ \blu{\calR'} := \{x\in \bbR^{nN} : u_i^\mpc (x_i,x_i^\tmp (P^*(x))) = -B_i^{-1} A_i x_i^\tmp (P^*(x)), \ \forall i\in \bbra{N} \} $.
Then, it holds
\begin{align}\label{eq:lyap_decrease_energy}
\frac{\rmd }{\rmd t} \calE (x(t), x^\sfd)
\begin{cases}
	< 0, & x(t) \not\in \blu{\calR'}, \\
	= 0, & x(t) \in \blu{\calR'}.
\end{cases}
\end{align}
In addition, by \eqref{eq:mpc_quad_conti}, we have
\begin{align*}
\blu{\calR'} &= \left\{ x\in \bbR^{nN} : x_i = x_i^\tmp (P^*(x)), \ \forall i \in \bbra{N} \right\} \\
&= \blu{\calR} .
\end{align*}
Finally, by applying \subm{again} the same argument as in the proof of Theorem~\ref{thm:lasalle}, we obtain the desired result.

\section{Proofs of Lemma~\ref{lemma:exponential_convergence} and Theorem~\ref{thm:stability}}\label{app:convergence}
\subsection{Proof of Lemma~\ref{lemma:exponential_convergence}}
Here, we abuse notation and identify $ P \in \bbR_{\ge 0}^{N\times N} $ as $ p = {\rm vec} (P) \in \bbR_{\ge 0}^{N^2} $.
It is known that the set of vertices of the Birkhoff polytope $ {\sf P} := \calT (\bfone_N/N, \bfone_N/N) $ is equal to the set of all permutation matrices~\cite{Birkhoff1946}.
Now, define
\begin{equation}\label{prob:LP_x}
	\xi (x) := \min_{P\in \calT(\frac{\one}{N}, \frac{\one}{N})} \sum_{i,j\in \bbra{N}} \bmC_{ij} (x) P_{ij} .
\end{equation}
Then, the set of optimal solutions of \eqref{prob:LP_x} is the intersection of $ {\sf P} $ and the hyperplane $ {\sf H} (x) := \{ P \in \bbR_{\ge 0}^{N\times N} : \sum_{i,j} \bmC_{ij} (x) P_{ij} = \xi(x) \}  $.

\blu{Note that since $ \bmA_i $ is invertible, $ \calG_i $ is also invertible.}
\blu{Then,} by $ \bmC_{ij}(x) = \|x_i - x_j^\sfd \|_{\calG_i}^2 $ and the assumption $ x_i^\sfd \neq x_j^\sfd $, for any $ \sigma $, the problem \eqref{prob:LP_x} with $ x = x^\sfd (\sigma) $ admits a unique optimal solution $ P^\sigma $, i.e., $ {\sf P} \cap {\sf H} (x^{\sfd} (\sigma)) = \{ P^\sigma \} $.
In addition, since the normal vector $ {\rm vec} (\bmC(x)) $ of the hyperplane $ {\sf H}(x) $ is continuous with respect to $ x $, we can take a neighborhood $ B_r (x^\sfd (\sigma)) $ where \eqref{prob:LP_x} with $ x \in B_r (x^\sfd(\sigma)) $ has the unique solution $ P^\sigma $.
\subm{By} the uniqueness and \subm{\cite[Proposition~5.1]{Cominetti1994}}, for any $ x \in B_r (x^\sfd (\sigma)) $, $ \bmP^* (x) $ converges exponentially to $ P^\sigma $ as $ \varepsilon \rightarrow +0 $.
Therefore, for any $ \delta > 0 $, we can choose sufficiently small $ \varepsilon > 0 $ such that
\begin{align}
	 &h(x) \in \calX_{\sigma, \delta}, \ \forall x \in B_r (x^\sfd (\sigma)), \label{eq:delta_cube} \\
	 &\calX_{\sigma,\delta} := \left\{x \in \bbR^{nN} : x^\sfd (\sigma) - \delta \bfone_{nN} \le x \le  x^\sfd (\sigma) + \delta \bfone_{nN} \right\} , \nonumber
\end{align}
where $ h $ is defined in \eqref{eq:func_fixed_point}, and the inequality sign between vectors should be understood element-wise.
Hence, \subm{by considering the restriction of $ h $ to $ \calX_{\sigma,\delta} $,} the same argument as in the proof of Proposition~\ref{prop:fixed_point_general} shows that for sufficiently small $ \delta $ and $ \varepsilon $, there exists at least one equilibrium $ x^\sfe (\varepsilon) \in \calX_{\sigma, \delta} \subset B_r (x^\sfd (\sigma)) $. \subm{In addition, $ x^\sfe (\varepsilon) $ converges to $ x^\sfd (\sigma) $ by letting $ \delta $ tend to zero, which implies $ \varepsilon \rightarrow + 0 $.}

Moreover, the exponential convergence of $ \bmP^* (x) $ for any $ x\in B_r (x^\sfd (\sigma)) $ implies that $ |\bmP_{ij}^* (x^\sfe (\varepsilon)) - P_{ij}^\sigma  | $ decays exponentially fast to $ 0 $ as $ \varepsilon \rightarrow + 0 $.
Lastly, for the convergence rate of $ x^\sfe (\varepsilon) $, we have
\begin{align*}
	\| x_i^\sfe (\varepsilon) - x_{\sigma(i)}^\sfd \| &= \biggl\| N\sum_{j=1}^N \bmP_{ij}^* (x^\sfe (\varepsilon)) x_j^\sfd - N \sum_{j=1}^N P_{ij}^\sigma x_j^\sfd \biggr\| \\
	&\le N\sum_{j=1}^N | \bmP_{ij}^* (x^\sfe (\varepsilon)) - P_{ij}^\sigma | \|x_j^\sfd \| ,  \forall i\in \bbra{N} .
\end{align*}
Thus, $ x^\sfe (\varepsilon) $ converges {\em exponentially} to $ x^\sfd (\sigma) $ as $ \varepsilon \rightarrow +0 $.

\subsection{Proof of Theorem~\ref{thm:stability}}
We prove only (ii) as the proof is similar for (i).
	In this proof, we regard $ (x (\cdot), \beta(\cdot) ) $ as a trajectory in a metric space $ \bbR^{nN} \times (\bbR_{>0}^N / {\sim} ) $ with \blu{the} metric $ d((x,\beta), (x',\beta')) := \|x - x'\| + \dhil (\beta, \beta') $.
	Fix any $ (x^\sfe, \beta^\sfe) \in {\rm Exp}(\sigma) $ satisfying the assumption in (ii).
	By definition, it is trivial that $ V $ \subm{in \eqref{eq:V}} is positive definite on a neighborhood of $ (x^\sfe, \beta^\sfe) $.
	Moreover, for any $ (x,\beta) \in \bbR^{nN} \times (\bbR_{>0}^N / {\sim})$, we have
	\begin{align*}
	&V(f_1 (x, \beta) , f_2 (x, \beta) ) - 	V(x , \beta ) \\
	&\le \sum_{i=1}^{N} \biggl\{ \Bigl\| \blu{\bar{\bmA}}_i \Bigl(x_i  -  N \sum_j \tilde{P}_{ij}(x,\beta)  x_j^{\sfd} \Bigr) \\
	&+ N\sum_j (\tilde{P}_{ij}(x,\beta) - P_{ ij}^\sfe  )x_j^{\sfd}  \Bigr\|_{\calG_i}^2 -  \Bigl\|x_i  - N\sum_j P_{ ij}^\sfe x_j^\sfd \Bigr\|_{\calG_i}^2 \biggr\}  \\
	&\quad + \gamma ( - W_3 (x,\beta) + W_4(x,\beta) + W_5 (x,\beta) ) \\
	&\le \sum_{i=1}^{N} \left(- W_{1,i}(x_i, \tilde{P}(x,\beta)) + W_{2,i}(x ,\beta) \right) \\
	&\qquad + \gamma ( - W_{3}(x,\beta) + W_{4}(x,\beta) + W_{5} (x,\beta)) =: W(x,\beta),
	\end{align*}
	where we used the triangle inequality for $ \dhil $, and
	\begin{align*}
		&W_{2,i} (x,\beta) := 2\Bigl( x_i  - N\sum_{j\in \bbra{N}} \tilde{P}_{ij}(x,\beta) x_j^\sfd \Bigr)^\top (\bar{\bmA}_i - I_n )^\top \calG_i  \\
		&\qquad\qquad\qquad \times N\sum_{j\in \bbra{N}} ( \tilde{P}_{ij}(x,\beta) - P_{ij}^\sfe) x_j^\sfd , \\
		&W_{3} (x,\beta) := \left[ 1- \lambda (\bmK(x)) \lambda \left(\bmK(f_1 (x,\beta)) \right) \right] \dhil(\beta, \beta^\sfe) , \\
		&W_{4} (x,\beta) :=  \dhil(\bmK(f_1 (x,\beta))^\top \alpha^\sfe, (K^\sfe)^\top \alpha^\sfe) , \\
		&W_{5} (x,\beta) := \lambda \left(\bmK(f_1(x,\beta)) \right) \dhil(\bmK(x)\beta^\sfe, K^\sfe \beta^\sfe), \\
		&K^\sfe := \bmK(x^\sfe), \ \alpha^\sfe := \bfone_N /N \oslash [K^\sfe \beta^\sfe].
	\end{align*}
	In the sequel, we explain that sufficiently small $ \varepsilon$ and large $\gamma $ enable us to take a neighborhood $ B_r(x^\sfe, \beta^\sfe) $ where
	\begin{equation}\label{ineq:W<0}
	  W(x,\beta) < 0, \ \forall (x,\beta) \in B_r(x^\sfe, \beta^\sfe)\backslash \{(x^\sfe, \beta^\sfe)\},
	\end{equation}
	which means the asymptotic stability of $ (x^\sfe, \beta^\sfe) $ \cite[Theorem~1.3]{Krabs2010}.

	First, a straightforward calculation yields, for any $ i,j\in \bbra{N}, l\in \bbra{n} $ and any $ (x,\beta) \in \bbR^{nN} \times (\bbR_{>0}^N / {\sim}) $,
	\begin{align*}
		&\left|\frac{\partial}{\partial x_{i,l}} \tilde{P}_{ij} (x, \beta) \right| \le \frac{2N \bar{g}_{i,j,l}}{\varepsilon} \tilde{P}_{ij} (x,\beta) \left( \frac{1}{N}  - \tilde{P}_{ij} (x,\beta) \right), \\
		&x_i = [x_{i,1} \ \cdots \ x_{i,n}]^\top, \\
		&\bar{g}_{i,j,l} := \max_{k\neq j} | g_{i,l}^\top (x_j^\sfd - x_k^\sfd) |, \ \calG_i = [g_{i,1} \ \cdots \ g_{i,n}]^\top .
	\end{align*}
	By Lemma~\ref{lemma:exponential_convergence}, under the assumption $ x_i^\sfd \neq x_j^\sfd, \ i\neq j $, $ \tilde{P}_{ij} (x^\sfe (\varepsilon), \beta^\sfe (\varepsilon)) $ converges exponentially to $ 0 $ or $ 1/N $ as $ \varepsilon \rightarrow +0 $.
	Hence, the variation of $ W_{2,i} $ with respect to $ x $ around $ (x^\sfe(\varepsilon), \beta^\sfe (\varepsilon)) $ can be made arbitrarily small by using sufficiently small $ \varepsilon = \varepsilon_1 $.
	 In addition, since $ \gamma > 0 $ can be chosen independently of $ \varepsilon $, sufficiently large $ \gamma = \bar{\gamma} $ enables us to take a neighborhood $ B_{r_1}(x^\sfe, \beta^\sfe) $ where
	 \begin{align}
		&\sum_{i=1}^N  \left(-\frac{1}{2}W_{1,i} \left(x_i, \tilde{P}(x,\beta) \right) + W_{2,i} (x,\beta)  \right) \nonumber\\
		& + \gamma \left(-\frac{1}{2}W_3(x,\beta) \right) < 0, \ 
	 	 \forall (x,\beta) \in B_{r_1}(x^\sfe, \beta^\sfe) \backslash \{(x^\sfe, \beta^\sfe)\}. \label{ineq:W_part1}
	 \end{align}

	Next, it follows from $ (x^\sfe, \beta^\sfe) \in {\rm Exp}(\sigma) $ that
	\begin{align*}
		\nabla_{x_i} \bmK_{ij}|_{x = x^{\las{\sfe}} (\varepsilon)} &= - \frac{2}{\varepsilon} \exp \left( - \frac{\|x_i^\sfe (\varepsilon) - x_j^\sfd \|_{\calG_i}^2}{\varepsilon} \right) \\
		&\qquad \times \calG_i (x_i^\sfe (\varepsilon) - x_j^\sfd) \rightarrow 0,  \ {\rm as} \ \varepsilon \rightarrow +0  .
	\end{align*}
	 Since $ W_4 $ and $ W_5 $ depend on $ (x,\beta) $ only via $ \bmK $, their variation around $ (x^\sfe (\varepsilon) ,\beta^\sfe (\varepsilon)) $ can be made arbitrarily small by taking sufficiently small $ \varepsilon > 0 $. Therefore, under the assumption that $ (x^\sfe(\varepsilon), \beta^\sfe (\varepsilon)) $ is isolated, for any given $ \gamma > 0 $, we can take $ \varepsilon = \varepsilon_2 (\gamma) $ such that there exists a neighborhood $ B_{r_2} (x^\sfe, \beta^\sfe) $ where
	 \begin{align}
		&\sum_{i=1}^N  \left(-\frac{1}{2}W_{1,i} (x_i, \tilde{P}(x,\beta))\right) \nonumber\\
		&\quad+ \gamma \biggl(-\frac{1}{2}W_3(x,\beta) + W_4(x,\beta)   + W_5(x,\beta) \biggr) < 0, \nonumber\\
	 	&\hspace{2cm} \forall (x,\beta) \in B_{r_2}(x^\sfe, \beta^\sfe) \backslash \{(x^\sfe, \beta^\sfe)\}. \label{ineq:W_part2}
	 \end{align}

	 By combining \eqref{ineq:W_part1} and \eqref{ineq:W_part2}, we obtain \eqref{ineq:W<0} for $ r = \min \{r_1, r_2\} $, $ \gamma = \bar{\gamma} $, and $\varepsilon = \min\{\varepsilon_1, \varepsilon_2 (\bar{\gamma})  \}$, which completes the proof.



\bibliographystyle{IEEEtran}
\bibliography{Sinkhorn_automatica}           

\end{document}